\documentclass[12pt]{amsart} 
\usepackage{color}
\usepackage{amssymb}
\usepackage{times}
\usepackage{amsmath}
\usepackage[pdftex]{graphicx}
\usepackage{hyperref}
\newcounter{sideremark}

\newtheorem{Thm}{Theorem}[section]
\newtheorem{Lem}[Thm]{Lemma}

\newenvironment{Pf}[1]
{\trivlist\item[]{\it #1\@. }}{\hspace*{\fill}$\Box$\endtrivlist}

\renewcommand{\marginpar}[1]{}
\catcode`\@=12

\def\Empty{}
\newcommand\oplabel[1]{
  \def\OpArg{#1} \ifx \OpArg\Empty {} \else
  	\label{#1}
  \fi}
		
\newcommand{\comm}[1]{}

\renewcommand{\epsilon}{\varepsilon}

\renewcommand{\rho}{\varrho}

\begin{document}
\title{Positive Ricci Curvature and the Length of a shortest periodic geodesic.}
\author{Regina Rotman}
\date{February, 23, 2022}
\maketitle

\begin{abstract}
Let $M^n$ be a closed Riemannian manifold of dimension  $n\geq 2$, 
with   Ricci curvature $Ric \geq n-1$. We will show
that any sphere of dimension $m$ 
in the space of closed loops on $M^n$ is homotopic
to the sphere in the space of closed loops of length at most 
$8 \pi m$. It follows that the length 
of a shortest periodic geodesic on $M^n$ is bounded from above 
by $8 \pi (n-1)$.
\end{abstract}


\section{Main results.}

Periodic geodesics are fundamental objects of Riemannian Geometry. In particular, it has been of long interest to bound the length of a shortest periodic 
geodesic in terms of the other geometric parameters of the manifold, such as its volume, diameter and $\backslash$  or curvature. Specifically, M. Gromov asked whether there exists a constant $c(n)$ such that on any closed Riemannian manifold of dimension $n$ and volume $vol$, the length of the shortest closed geodesic can be majorized by 
$c(n) vol^{\frac{1}{n}}$, (see [G]). A similar question can be posed in terms of  diameter of the manifold instead of its volume. The above questions were well-studied in dimension $2$, beginning with  results of K. Loewner and P. Pu, who have established sharp upper bounds in cases of a Riemannian $2$-torus and a Riemannian $RP^2$ respectively. Some other  notable contributions include those of J. Hebda, Yu. Burago and V. Zalgaller, A. Treibergs, M. Gromov, C. Bavard, (who proved a sharp upper bound in the case of Riemannian Klein bottle), as well as M. Katz and S. Sabourau, (see [Ba], [H], [G],  [KS1], [KS2], as well as [BrZ], [CrK], [K] for  surveys of some of the results). Among the surfaces, the most challenging case was 
that  of a Riemannian $2$-sphere. The first area and diameter upper bounds in this case were found 
by C. B. Croke in [Cr], where it was shown that in the case of a Riemannian $2$-sphere $M$ of area $A$ and diameter $d$, the length of the shortest periodic geodesic $l(M)$ can be majorized by $31 \sqrt{A}$ and $9d$. Those bounds were later improved in [M], [S], [NR1], and [R1]. A local optimal area upper bound for the length of a closed periodic geodesic on the 2-sphere in a neighbourhood of the singular metric obtained by gluing two equilateral triangles along their boundary was established in [B1]. 
The only known
curvature-free estimate in dimension higher than two is the result of 
M. Gromov for the class of $1$-essential manifolds (see [G] and also [N], which contains the best currently known bound). 

However, there were established a  few upper
bounds for the length of a shortest periodic geodesic in terms of 
sectional and Ricci curvature bounds, starting from the results of W. Ballman, G. Thorbergsson, W. Ziller, ([BlTZl]), who dealt with the case of $\frac{1}{4}$-pinched metric of positive sectional curvature on spheres.
Other prior results include the case of convex hypersurfaces, due to Croke, ([Cr]) and Treibergs, ([T]), as well as
the more recent result for the metrics on closed Riemannian manifolds with diameter bounded above, volume 
 bounded below and sectional curvature bounded below, (see [NR2]), and the very recent results of N. Wu and Z. Zhu
dealing with the case of manifolds of dimension four with the bounded Ricci curvature [WZ], as well as the result of I. Adelstein and F. Vargas Pallete ([AVP]), in which they further improve the diameter upper bound for the length of the shortest closed geodesic on a positively curved Riemannian $2$-sphere.

In this paper we will prove the following theorem. 

\begin{Thm} \label{TheoremAA}
Let $M^n$ be a closed Riemannian manifold of dimension $n$ and 
$Ric \geq n-1$. Let $f: S^m \longrightarrow \Omega_p M^n$ be (a map of) a
sphere of dimension $m \geq 1$ in the space of all piecewise differentiable loops based at some point $p \in M^n$, $\Omega_p M^n$ on $M^n$. Then there exists a homotopy 
$H: S^m \times [0,1]$ such that $H(*,0)=f$ and $H(*,1) = \tilde{f}$, 
where $\tilde{f}: S^m \longrightarrow \Omega_p^{8 \pi m}M^n$. Here
$\Omega_p^{8 \pi m}M^n$ is the space of piecewise differentiable loops based at $p$ of length at most $8 \pi m$ on $M^n$.  Moreover, if the image of the original map $f$ is contained inside $\Omega_p^LM^n$, one 
can choose $H$ so that its image is contained in $\Omega_p^{L+8\pi m+o(1)}$.
\end{Thm}

\noindent{\bf Remark:} It is easy to see that for closed Riemannian 
manifolds with $Ric \geq (n-1)$ and $\pi_1(M^n) \neq \{0\}$, the length of 
a shortest periodic geodesic, $l(M^n) \leq \pi$. In fact, let {\it C} be a free homotopy class. Let $\alpha=\inf length (\gamma), \gamma \in $ {\it C}. (The proof of) Bonnet-Myers Theorem implies that the index of any piecewise differentiable
curve of length $> \pi$ is at least $1$. In particular, any piecewise differentiable closed curve of length $>\pi $ can be continuously shortened. 
Therefore, $\alpha \leq \pi$. By Cartan's theorem, there exists a closed 
curve $\tilde{\gamma} \in $ {\it C} of length $\alpha$, which is a closed 
geodesic, (see [DC]).

In 1959 it was proven by A. Fet and L. Lusternik that on any closed
Riemannian manifold $M^n$ there exists at least one periodic geodesic. 
The proof uses Morse theory on the space $\Lambda M^n$  of piecewise differentiable closed curves on $M^n$. From the homotopy theory it follows that 
there exists an integer $q$ such that $\pi_q(M^n) = \{0\}$, while 
$\pi_q(\Lambda M^n) \neq \{0\}$. Note that we can identify points on $M^n$ 
with constant curves on $M^n$. Let us consider the Energy functional, $E$ 
on $\Lambda M^n$, and let $f: S^q \longrightarrow \Lambda M^n$ be a 
non-contractible sphere in $\Lambda M^n$. Let us try to deform this sphere 
to $M^n$, identified with the space of constant curves on $M^n$,  along the integral curves of $\nabla E$. This should be impossible, since if one succeeds, one can 
then contract this sphere in $M^n$, thus reaching a contradiction. Therefore, 
there exists a critical point of $E$ that should prevent one from contracting 
this sphere to $M^n$. This critical point is a periodic geodesic on $M^n$. 

We can see that in conjunction with Theorem ~\ref{TheoremAA} the proof 
of Fet and Lusternik theorem implies the following result: 

\begin{Thm} \label{Theoremmain}
Let $M^n$ be a closed Riemannian manifold of dimension $n$ with the 
Ricci curvature $Ric \geq n-1$. Let $q$ be the smallest integer such that
$\pi_{q+1}(M^n) \neq \{0\}$. 

Then there exists a periodic geodesic 
on $M^n$ of length at most $8 \pi q \leq 8 \pi (n-1)$. 
\end{Thm}

From the above result it immediately follows that the length of 
a shortest periodic geodesic on any Riemannian manifold 
with Ricci curvature bounded from below by $\frac{n-1}{r^2}$ for some
positive $r$ is bounded from above by $8\pi r q$, where 
$q$  is the smallest integer such that $\pi_q(M^n) =\{0\}$, 
$\pi_{q+1}(M^n) \neq \{0\}$. 

Another set of related questions deals with bounding the length of such  objects as critical points of the length functional on the space of graphs on a manifold, on the space of closed curves based at a given point and on the space of curves between a fixed pair of points. Those objects are geodesic nets, geodesic loops and geodesics
respectively. Many curvature-free upper bounds are known in the case of nets, loops and segments, starting from the 
result of S. Sabourau, who established the first curvature-free diameter and volume upper bounds for the length of the shortest geodesic loop on a closed Riemannian manifold, (see [S2]). For other results of similar nature
see [B], [NR3], [NR4], [NR5], [NR6], [NR7], [R2], [R3], [R4], [R5]. 
Another result that immediately follows from Theorem ~\ref{TheoremAA} is the following quantitative version of J. P. Serre's theorem that states that for any 
pair of points $p, q \in M^n$ of a closed Riemannian manifold $M^n$ there exist infinitely many geodesics that connect these points. Note that when $p=q$, this result implies that at any point $p \in M^n$, there exists infinitely many geodesic loops, (see [Se]).  The first effective proof of this theorem was presented by 
A. Schwartz, who in particular demonstrated that given any pair of points $p, q \in M^n$, where $M^n$ is a closed Riemannian manifold, there exists at least $k$ geodesics connecting them of length at most $c(M^n)k$, where 
$c(M^n)$ depends on a Riemannian metric on $M^n$, (see [Sc]). Together with A. Nabutovsky, we have shown 
that in the case of a Riemannian $2$-sphere, one can take $c(M^n)=22$, ($20$ when $p=q)$, (see [NR5], [NR6]). These bounds were recently  significantly improved by H. Y. Cheng in [Ch].  

\begin{Thm} \label{Theoremmain2}
Let $M^n$ be a closed Riemannian manifold of dimension $n$ with Ricci curvature $\geq n-1$. Then 

\noindent (1) given any point $p \in M^n$
there exists at least $k$ geodesic loops based at $p$ of length at most 
$16\pi (n-1) k$;

\noindent (2) given a pair of points $p, q \in M^n$ there exist at least $k$ geodesics of length at most $\pi (16k(n-1)+1)$.
\end{Thm}

Theorem ~\ref{Theoremmain2} improves the previously known bound for the length of the first
$k$ geodesic loops based at  $p$ from [NR7] in the  case of the manifolds with a positive Ricci curvature. While in the theorem in [NR7] that  states that on any closed Riemannian manifold there exist at least $k$ geodesic loops of length at most $4nk^2d$, the bound is quadratic in $k$, in case of the manifolds with positive Ricci curvature, one can improve this to be a linear bound in $k$.  


\section{Main Ideas}

Let $M^n$ be a closed Riemannian manifold. The existence of a closed geodesic theorem of Fet and Lyusternik, as well as that of infinitely many geodesics between a pair of points of Serre rely on the Morse Theory on $\Lambda M^n$ and $\Omega_{p,q}M^n$ respectively. Thus, as we have demonstrated the quantitative versions of these theorems can be obtained by constructing of the short spheres of an 
appropriate dimension in the given spaces of curves. This, in turn is obtained from obtaining a relevant sweep-out of the manifold. While, it has been demonstrated that one cannot hope to construct such sweep-outs in the general situation, (see [FK], [L]), usually such sweep-outs are prevented by the 
existence of  stationary objects of index zero.  However, in our case of manifolds with a positive 
Ricci curvature, we demonstrate such such short sweep-outs are indeed possible. 

The proof of Theorem ~\ref{TheoremAA} has the following two main ingredients. 

\noindent (1)  The first ingredient is Bonnet- Myers Theorem, 
which implies that when  $M^n$ is a Riemannian manifold 
with $Ric \geq n-1$, the diameter $d$ of $M^n$ is $\leq \pi$. 
It also implies that geodesics of length $>\pi$ between any pair of 
points $p, q \in M^n$ has index $\geq 1$. This, in particular, is 
true for geodesic loops.

\noindent (2) The second ingredient is the idea introduced in 
[NR7] that can be informally phrased as the following observation:
whenever one can shorten fixed point loops via loops of controlled length, one can also shorten families of 
fixed point loops over loops of controlled length. We would like to stress, that the families of loops 
can be shortened even if the fixed point  homotopies are not  continuous with respect to the original curves. In other words, it can happen that two ``close'' loops are shortened by homotopies that are ``far'' from each other. 

The following theorem (Theorem ~\ref{TheoremAAA}) is a slight modification of Theorem 1.1 in [NR7], that 
can be restated to say that if for some positive integer $k$ and some point 
$p$ in a closed Riemannian manifold $M^n$ no locally minimizing geodesic loops
based at $p$ have length in the interval $(2(k-1)d,2kd]$, then every map 
$f:S^m \longrightarrow \Omega_pM^n$ is homotopic to $\tilde{f}:S^m \longrightarrow \Omega_p^LM^n$, where $L=((4k+2)m+(2k-3))d+\upsilon$. Here
$d$ is the diameter of $M^n$, $n$ is the dimension, $\Omega_pM^n$ denotes the space of piecewise differentiable loops  based at $p$, whereas $\Omega_p^LM^n$ denotes the of piecewise differentiable loops based at $p$ of length at most $L$.   While in the statement of the above theorem, we assumed that $k$ is a positive integer, it is easy to see that the conclusion of the theorem will also 
hold for any positive real number $k$.  

\noindent {\bf Remark.} Note that if we then take $k=\frac{\pi}{2d}+1$ ($d$ 
being the diameter of $M^n$ as before) then by Bonnet-Myers Theorem, 
there will be no locally minimizing loops based at any $p \in M^n$ with 
length greater than $\pi$. Thus, this relaxed hypothesis of Theorem 1.1 in 
[NR2] is satisfied for such a $k$. It then immediately follows 
that for an arbitrary small positive $\upsilon$ the length of a shortest periodic geodesic on $M^n$ is at most 

$[(4(\frac{\pi}{2d}+1)+2)q+(2(\frac{\pi}{2d}+1)-3)]d \leq
[({2\pi \over d}+6)q+({\pi \over d}-1)]d + \upsilon \leq
8 \pi q+\pi - d + \upsilon \leq 9 \pi q + \upsilon$. This is only a slightly worse bound then stipulated in our theorem.   
       
One can, however, use methods of [NR7] to prove a slightly stronger results than Theorem 1.1 in [NR7].

\begin{Thm} \label{TheoremAAA} Let $M^n$ be a closed Riemannian manifold of 
dimension $n$ and diameter $d$. Suppose $d \leq a$ for some positive
real number $a$. Let $p \in M^n$ be an arbitrary point of $M^n$ and $k$ be
some positive real number 
 for which no geodesic loop based at $p$ of index zero
has length $l$ such that $2(k-1)a < l \leq 2ka$. Then every map 
$f: S^m \longrightarrow \Omega_pM^n$ is homotopic to a map 
$\tilde{f}:S^m \longrightarrow \Omega_p^LM^n$, where 
$L=((4k+2)m+(2k-3))a+\upsilon$, where $\upsilon$ is an arbitrarily small 
positive number. 
\end{Thm}

Then Theorems \ref{TheoremAA} and \ref{Theoremmain} immediately follow from 
Theorem \ref{TheoremAAA}:

\begin{Pf}{Proof of Theorems \ref{TheoremAA} and \ref{Theoremmain}}  As we observed, Bonnet-Myers theorem
implies that closed Riemannian manifolds
with $Ric \geq n-1$ satisfy the hypothesis of 
Theorem ~\ref{TheoremAAA} with $k={3 \over 2}$ and $a=\pi$. Thus, we obtain Theorem
\ref{TheoremAA}. Now let $q$ be the smallest positive integer number such that $\pi_{q+1}(M^n)$ is non-trivial. Consider a non-contractible map $S^{q+1}\longrightarrow M^n$
and the corresponding map $S^q\longrightarrow \Omega_pM^n$, where $p$ is any point in $M^n$. We saw that if $q=0$, the length of the shortest periodic geodesic
does not exceed $\pi$, so we can assume that $q\geq 1$.
Apply Theorem \ref{TheoremAA} to conclude that this map is homotopic to a map
$S^q\longrightarrow \Omega_p^{8\pi q}M^n\subset \Lambda^{8\pi q}M^n$. Now the proof of the Fet-Lyusternik theorem immediately implies that
the length of a shortest periodic geodesic on $M^n$ is 
at most $8 \pi q$. 
\end{Pf}

Now we will sketch the proof of Theorem \ref{Theoremmain2} from Theorem \ref{TheoremAAA}.

\begin{Pf}{Proof of Theorem \ref{Theoremmain2}}The proof of Theorem ~\ref{Theoremmain2} likewise directly follows from Theorem ~\ref{TheoremAA} as follows:
Rational homotopy theory implies the existence of a non-trivial 
even-dimensional cohomology class of dimension $\leq 2(n-1)$ on a path/ loop space on a closed Riemannian manifold $M^n$, whose cup powers are not trivial, (see [NR7], [FHT]). The classes generated by those cup powers will correspond to the different geodesic paths / loops in a Morse situation. Otherwise, if the two classes happen to correspond to the same critical point, Lyusternik and Schnirellmann theory would imply the existence of the whole critical level. As it was further demonstrated by A. Schwartz in 
[Sc], one can consider dual homology classes under the Pontryagin powers instead of the cohomology classes under the cup powers. Thus, if one can represent the corresponding class by short segments / loops, it will further imply the existence of short paths / loops in all dimensions, (see [NR7]). Now we again observe that
Bonnet-Myers theorem and Theorem \ref{TheoremAAA} imply that this class can be represented by loops of length $\leq 16\pi (n-1)$ (correspondingly, paths $\leq \pi (16(n-1)+1)$). The Pontryagin powers can be represented by loops of length $\leq 16\pi(n-1)k$ (correspondingly, paths of length $\leq \pi(16(n-1)k+1)$).
\end{Pf}


\section{Summary of the proof of Theorem ~\ref{TheoremAAA}}

The proof of Theorem ~\ref{TheoremAAA} is very similar to the proof of 
Theorem 1.1 of [NR7]. Therefore, only the summary of the proof will 
be included here for the paper to be self-contained.

The starting point of the proof is the following effective version of 
the elementary topological lemma (Lemma 2.1 from [NR7]) and its higher dimensional generalization (Lemma 2.2 from [NR7]). Let $\alpha:[0,1] \longrightarrow M^n$ be a piecewise differentiable path in $M^n$. Then $\bar{\alpha}$ will denote
$\alpha(1-t)$ in the statement of the lemma and thereafter in the paper. 
Moreover, if $\beta$ is another path in $M^n$, such that $\alpha(1)=\beta(0)$, 
then $\alpha*\beta$ will denote their product. 

Informally, the following lemma can be viewed as follows: consider two edges 
starting and ending at the same pair of points $p, q$. By looking at a concatenation of these edges we can view them as a loop based at $p$. Then if this loop is path homotopic to another loop based at $p$ via loops of short length, there is a ``short'' path homotopy between the first edge and a concatenation of the  loop followed by the second edge. 

\begin{Lem} \label{Lemma}   
Let $\gamma_1$, $\gamma_2$ be two paths on a complete Riemannian manifold $M^n$  of length $l_1, l_2$ respectively, such that $\gamma_1(0)=\gamma_2(0)=p$, 
and $\gamma_1(1)=\gamma_2(1)=q$. 
Let $\alpha=\alpha_0=\gamma_1 * \bar{\gamma}_2$. This is a loop based 
at point  $p$. Suppose further that there exists a path homotopy   
$H(t,\tau)$ between $\alpha_0$ and a loop $\alpha_1$ over 
the loops $\alpha_\tau$ of length at most $l_3$. Then 
there exists a path homotopy $\tilde{H}(t,\tau)$ between $\gamma_1$ and 
$\alpha_1 * \gamma_2$ over curves of length at most $l_3+l_2$. Similarly, 
$\gamma_2$ is path homotopic to $\bar{\alpha}_1 * \gamma_1$ over the curves 
of length at most $l_1+l_3$. 
\end{Lem}

\begin{figure}[!htbp]
\centering
\includegraphics[width=9cm]{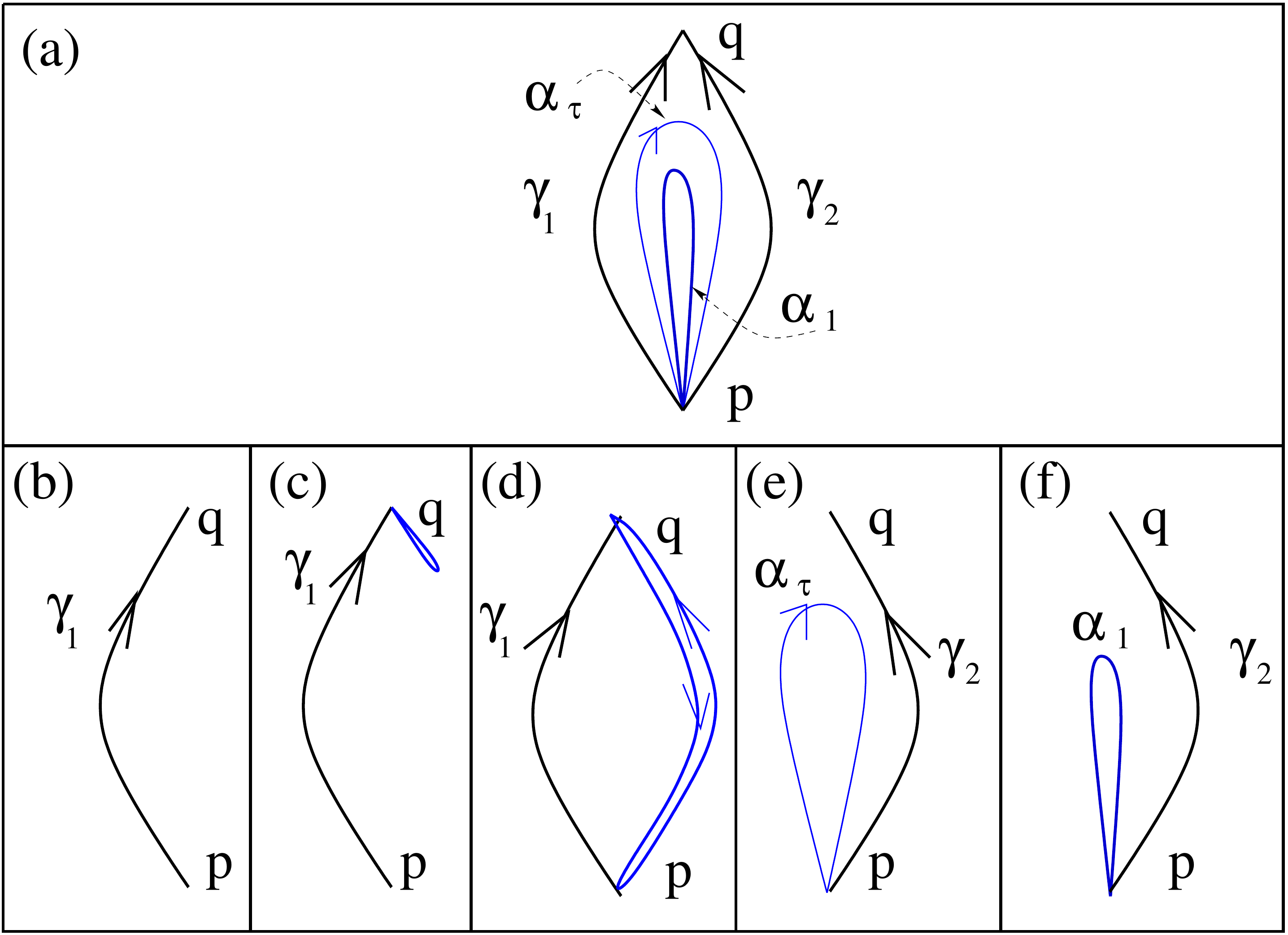}
\caption{Illustration of the proof of Lemma~\ref{Lemma}}
\label{figure1}
\end{figure}

\begin{Pf}{Proof}
We are going to construct a homotopy $\tilde{H}(t, \tau)$ between  
$\gamma_1$ and $\alpha_1 * \gamma_2$ as follows. 
For the illustration of the proof see Fig. ~\ref{figure1}.  
First note that 
$\gamma_1$ is path homotopic to $\gamma_1*\bar{\gamma}_2*\gamma_2$. The typical curve in this homotopy consists of $\gamma_1$ followed by a  longer and longer segment of $\gamma_2$ taken twice with the opposite orientation. That is,  one goes a short distance from $p$ towards $q$ along $\gamma_2$, and then goes back  to $q$, (see Fig.  ~\ref{figure1} (b), (c), (d)). The  curves in this homotopy are bounded by $l_1+2l_2 \leq l_3+l_2$.  Next, 
$\gamma_1*\bar{\gamma}_2$ is path homotopic to $\alpha_1$ via the curves $\alpha_\tau$
of length $\leq l_3$. Therefore,  the path
$\gamma_1*\bar{\gamma}_2*\gamma_2$ is path homotopic to $\alpha_1*\gamma_2$ along the curves
$\alpha_\tau*\gamma_2$ of length at most $l_2+l_3$; (see Fig. ~\ref{figure1} (d,e)).
\end{Pf}

{\bf Remark:} One can apply Lemma ~\ref{Lemma} in the situation, when $\alpha_1$ is a constant path at $p$. In this case  the resulting homotopy will be a short path homotopy between $\gamma_1$ and $\gamma_2$. 

Below we will state the higher dimensional generalization of this lemma.

Let $p, q \in M^n$ be a  pair of points on a complete Riemmanian manifold of dimension $n$, and  let $\Omega_{p,q}^LM^n$ be  the space of piecewise differentiable curves
on $M^n$ of length at most $L$.  

\begin{figure}[! htbp]
\centering\includegraphics[width=9cm]{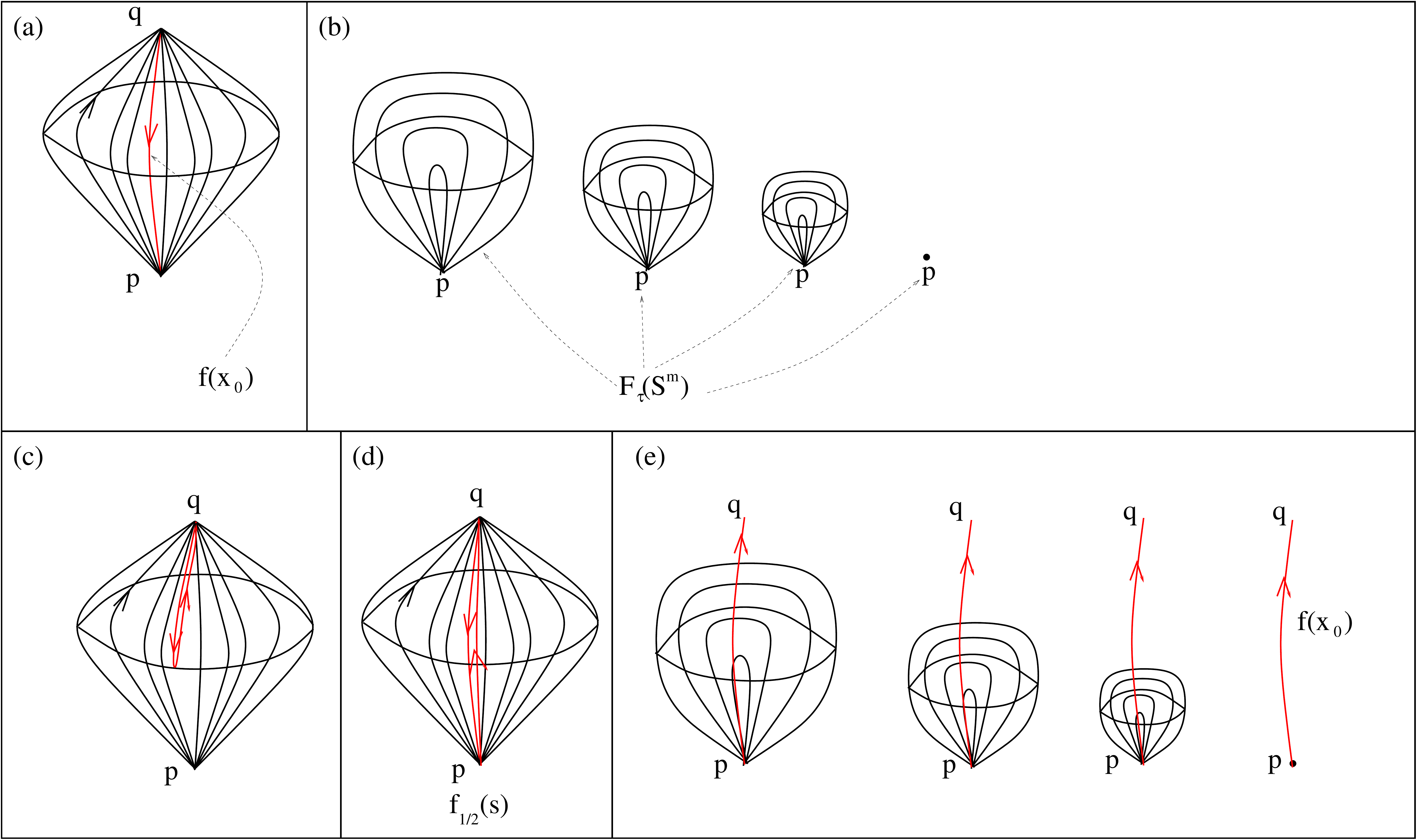}
\caption{Illustration of the proof of Lemma~\ref{Lemmagen}}
\label{figure2}
\end{figure}


\begin{Lem} \label{Lemmagen}
Let $f: S^m \longrightarrow \Omega^L_{p,q}M^n$ be a continuous map from the 
m-dimensional sphere into $\Omega^L_{p,q}M^n$. Let $l$ be  the length of
the curve  $f(x_0)$ for some $x_0 \in S^m$ (see Fig. ~\ref{figure2} (a)).


Consider a sphere in the space of loops based at $p$ that is 
formed by concatenating the paths in the image of $f$ by  $\overline{f(x_0)}$. That is we define 
 $F:S^m \longrightarrow
\Omega_p^{L+l}M^n$ as $F(x)=f(x)*\overline{f(x_0)}$. Then, if there exists 
a homotopy $F_\tau:S^m \longrightarrow \Omega_p^{L+l}$ contracting the sphere
to $p$  (See Fig. ~\ref{figure2} (b)),  then
there exists a homotopy $f_\tau:S^m\longrightarrow
\Omega_{p,q}^{L+2l}M^n$, $\tau \in [0,1]$,
between $f=f_0$ and the constant map $f_1$ of $S^m$
identically equal to $f(x_0)$ (see Fig. ~\ref{figure2} (c)-(e)).
\end{Lem}

\begin{Pf}{Proof} 
The homotopy is constructed first by attaching to each curve $f(x)$ longer and longer segment of $f(x_0)$ taken twice with different orientations. Eventually, we will attach the curve 
$\overline{f_(x_0)} * f(x_0)$. 
Formally speaking,   $f$ is homotopic to  $f_{1\over 2}$ defined as
$f_{1\over 2}(s)=f(x)*\overline{}{f(x_0)}*f(x_0)=F(x)*f(x_0)$. 
\par
Next we contract the sphere $F$ using the homotopy $F_t,\ t\in [0,1]$, while the curve $f(x_0)$
that is attached to each loop stays fixed, (see Fig. ~\ref{figure2} (c)).
\end{Pf}

\begin{figure}[! htbp]
\centering\includegraphics[width=9cm]{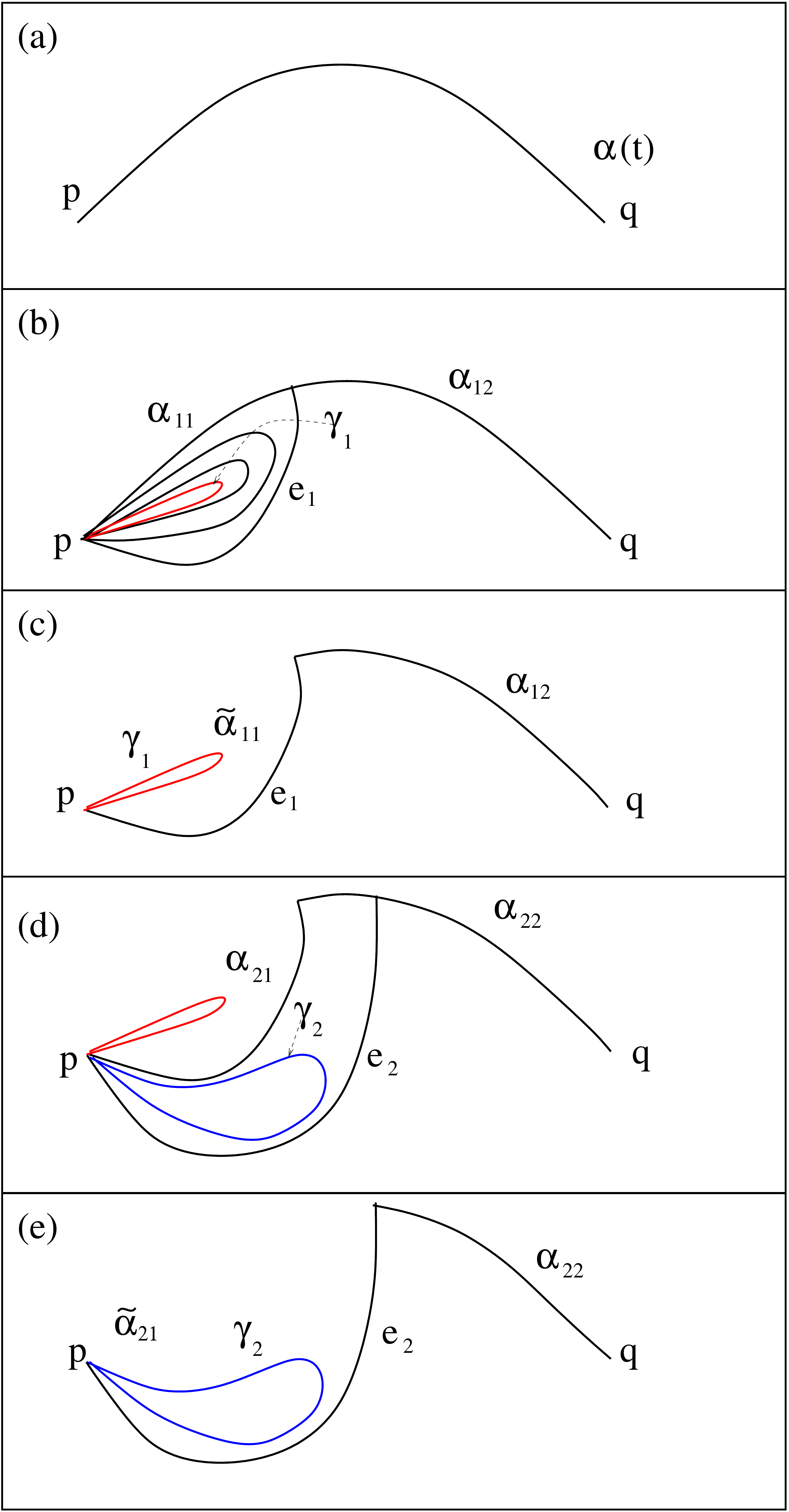}
\caption{Curve shortening process}
\label{figure3}
\end{figure}

The main tool used in  Theorem ~\ref{TheoremAA} is  a 
curve shortening process described in the following theorem

\begin{Thm} \label{TheoremA}
Let $M^n$ be a closed Riemannian manifold of
diameter $d$. Suppose that $d \leq a$ for some positive real number 
$a$. Let $\alpha = \alpha(t)$ be a curve on $M^n$ connecting an arbitrary pair of points
$p, q \in M^n$. Assume that
there exists an interval $(l, l+2a+\delta]$, such that there are no
geodesic loops based at $p \in M^n$  with length in this interval of index zero with respect to the length functional on $\Omega_pM^n$ for some positive $\delta$. Then there 
exists a curve $\tilde{\alpha} = \tilde{\alpha}(t)$ of length $\leq l+a$ 
connecting $p$ and $q$ and a path homotopy between $\alpha$
and $\tilde{\alpha}$ such that
the lengths of all curves in this path homotopy are at most $L+2a$.

Moreover, there exists a one-parameter family of curves: $\gamma_s^\alpha$
that continuously connect $p$ with points of $\alpha(t)$ satisfying the following conditions:

\noindent (a) length $(\gamma_s^\alpha) \leq l+3a+\delta$;

\noindent (b) for every $s$, $\gamma_s^\alpha$ connects $p$ with
$\alpha(\tau(s))$, where $\tau(s)$ is a function of $s$ satisfying
$\tau(s_1) \leq \tau(s_2)$ whenever $s_1 < s_2$;

\noindent (c) There exists two partitions of $[0,1]$,
$P=\{0=t_0^\alpha < t_1^\alpha< ... < t_{k^\alpha}^\alpha=1\}$ and
$Q=\{0=r_0^\alpha < r_1^\alpha < ... <r_{{2k}^\alpha}^\alpha=1\}$, such that

\noindent (1) the endpoint of $\gamma_s^\alpha$ remains constant for
$s \in [r_{2i-1}^\alpha, r_{2i}^\alpha]$;

\noindent (2) for every $s \in [r_{2i}^\alpha, r_{2i+1}^\alpha]$, the endpoint
 of $\gamma_s^\alpha$ is $\alpha(t)$ for some $t \in [t_i^\alpha, t_{i+1}^\alpha]$.  Moreover, $\tau(s)$ is strictly increasing on $[r_{2i}^\alpha, r_{2i+1}^\alpha]$, and $\tau([r_{2i}, r_{2i+1}])=[t_i,t_{i+1}]$.

\noindent (3)  the length of $\gamma^\alpha_{s^\alpha_{2i}} \leq l+a$ for all \
$i$.

\noindent (4) the curve $\gamma^\alpha_{r_{2i}}=\gamma_1^\alpha$ is the final curve in the curve shortening above.

\end{Thm}

\begin{figure}[! htbp]
\centering
\includegraphics[width=9cm]{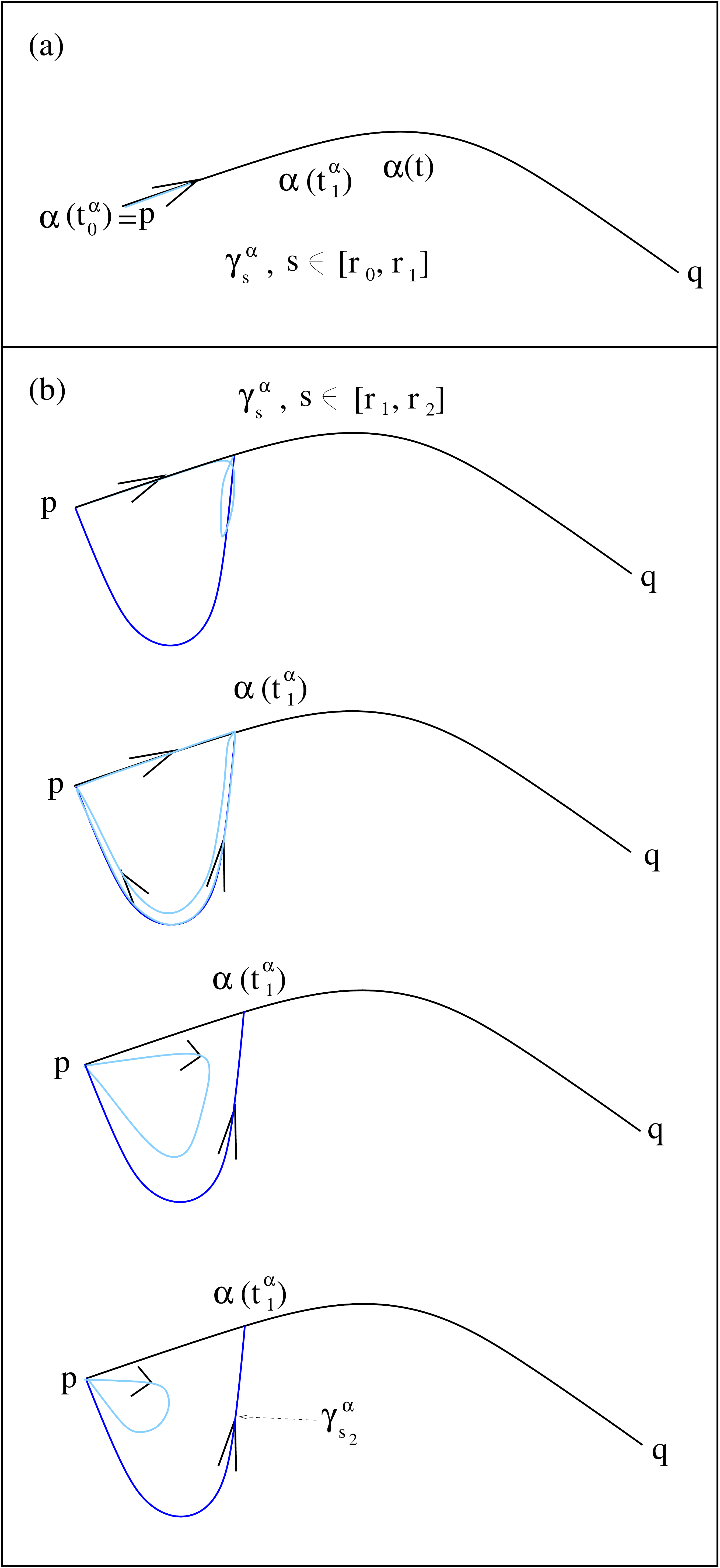}
\caption{Constructing $1$-parameter family of curves: $\gamma^\alpha_s$}
\label{figure4}
\end{figure}

\par

\begin{Pf} {Proof of Theorem ~\ref{TheoremA}}

Let $\alpha:[0,L] \longrightarrow M^n$ be a curve of length $L$ parametrized by the arclength (see Fig. ~\ref{figure3} (a)). If $L \leq l+a$, then we are done. Suppose  $L > l+a$. 
Consider the segment $\alpha\vert_{[0,l+a+\delta]}$, where $\delta$ is specified in the hypothesis of the theorem, (or possibly smaller, if the whole curve $\alpha$ is shorter than that) of $\alpha$,
which we will denote $\alpha_{11}=\alpha_{11}(t)$, and the segment $\alpha\vert_{[l+a+\delta, L]}$
denoted $\alpha_{12}=\alpha_{12}(t)$
(see Fig. ~\ref{figure3} (b)).

Let us connect points $p$ and $\alpha(l+a+\delta)$ by a minimal 
geodesic, $e_1$ (of length $\leq d \leq a$) (see Fig. ~\ref{figure3} (b)). 
Then curves
$\alpha_{11}$ and $e_1$ form a loop $\alpha_{11}*\bar e_1$
of length $\leq l+2a+\delta$
based at $p$.

Consider a (possibly trivial) shortest loop $\gamma_1$ of index zero that
can be connected with $\alpha_{11}*\bar e_1$ by a length non-increasing
path homotopy, (see Fig. ~\ref{figure3} (b)).  The length of  $\gamma_1$ is at most $l$ by our assumption.

By Lemma ~\ref{Lemma}, $\alpha_{11}$ is
path homotopic to the curve $\gamma_1*e_1=\tilde{\alpha}_{11}$
along the curves of length at most $l+3a+\delta$ and, thus,
the original curve $\alpha$ is homotopic to a new curve $\tilde{\alpha}_{11}*\alpha_{12}$ along the curves of length at most $L+2a$, (see Fig. ~\ref{figure3} (c)). 

Note that the 
length $L_1$ of this new curve
$\alpha_1=\tilde{\alpha}_{11}*\alpha_{12}$
is at most $L -\delta$.  Assuming that $L_1$
is still greater than $l+a$, we repeat the process again:  
We parametrize $\alpha_1$ by its arclength.
Now, let $\alpha_{21}=\alpha_1\vert_{[0,l+d+\delta]}$
and $\alpha_{22}=\alpha_1\vert_{[l+a+\delta, L_1]}$, (see Fig. ~\ref{figure3} (d)). If it so happens that the length $L_1$ of $\alpha_1$ is less than $l+a+\delta$, 
then
we use $L_1-l-a$ as the new value of $\delta$, but otherwise
keep the old value of $\delta$.

Connect the points
$p$ and $\alpha_1(l+d+\delta)$ by a minimal geodesic segment $e_2$,
(see Fig. ~\ref{figure3} (d)).
Then $\alpha_{21}$ and $e_2$ form a  loop $\alpha_{21}*\bar e_2$
based at $p$
of length at most $l+2a+\delta$.  

Again, we use the length shortening to contract this loop to the 
shortest possible loop $\gamma_2$ of index $0$ keeping the vertex fixed, (see Fig. ~\ref{figure3} (e)).
By our assumption the length of $\gamma_2$ is at most $l$.

Thus, $\alpha_{21}$ is path-homotopic to 
$\tilde{\alpha}_{21}=\gamma_2*\sigma_2$
along the curves of length at most $l+3a+\delta$.  It follows
that $\alpha_1$ is homotopic to $\alpha_2=\tilde{\alpha}_{21}*\alpha_{22}$
along the curves of length at most $L+2a$, (see Fig. ~\ref{figure3} (e)).

Let $k$ be the number of steps to sufficiently shorten $\alpha(t)$. One can see that $k= \lfloor{\frac{L-l-a}{\delta}}\rfloor+1$. 

In the process of shortening the original curve $\alpha(t)$ we are constructing the axillary curves $\gamma_s^\alpha$ in the following way: we connect $p$ with the points of $\alpha(t)$ by minimizing geodesics: $e_1, e_2,..., e_k$.  Each points that are distance $\delta$ apart on $\alpha(t)$ are being connected, started from the point $\alpha(l+a+\delta)$. We next construct $\gamma_s^\alpha$ by repeated application of 
Lemma ~\ref{Lemma}. As we have shown in the process of the proof, the  resulting family of curves $\gamma_s^\alpha$
continuously connects $p$ with points of $\alpha(t)$ satisfying the following conditions:

\noindent (a) length $(\gamma_s^\alpha) \leq l+3a+\delta$;

\noindent (b) for every $s$, $\gamma_s^\alpha$ connects $p$ with 
$\alpha(\tau(s))$, where $\tau(s)$ is a function of $s$ satisfying 
$\tau(s_1) \leq \tau(s_2)$ whenever $s_1 < s_2$;

\noindent (c) There exists two partitions of $[0,1]$, 
$P=\{0=t_0^\alpha < t_1^\alpha < ... < t_{k^\alpha}^\alpha=1\}$ and 
$Q=\{0=r_0^\alpha < r_1^\alpha < ... <r_{{2k}^\alpha}^{\alpha}=1 \}$, such that

\noindent (1) the endpoint of $\gamma_s^\alpha$ remains constant for 
$s \in [r_{2i-1}^\alpha, r_{2i}^\alpha]$;

\noindent (2) for every $s \in [r_{2i}^\alpha, r_{2i+1}^\alpha]$, the endpoint of $\gamma_s^\alpha$ is $\alpha(t)$ for some $t \in [t_i^\alpha, t_{i+1}^\alpha]$.  Moreover, $\tau(s)$ is strictly increasing on $[r_{2i}^\alpha, r_{2i+1}^\alpha]$, and $\tau([r_{2i}, r_{2i+1}])=[t_i,t_{i+1}]$. 

\noindent (3)  the length of $\gamma^\alpha_{s^\alpha_{2i}} \leq l+a$ for all $i$, since each of those curves consist of a geodesic $e_i$ of length at most $a$  with a loop attached of length at most $l$. 

\noindent (4) the curve $\gamma^\alpha_{r_2i}=\gamma_1^\alpha$ is the final curve in the curve shortening above.

\end{Pf}

\par

On Figure ~\ref{figure4}, we have depicted the curves $\gamma^\alpha_s$, where
$s \in [r_0, r_1]$, (Fig. ~\ref{figure4}(a)) and for $s \in [r_1, r_2]$, (Fig. ~\ref{figure4} (b)). Note that for $s \in [r_0, r_1]$, $\gamma_s^\alpha$ will be just longer and longer segments of $\alpha$, whereas $\gamma_s^\alpha, s \in [r_1, r_2]$ will be obtained via
Lemma ~\ref{Lemma}. Therefore, they will be the curves interpolating between 
$\alpha_{11}$ and $\gamma_1 \star e_1$. Note also that $\gamma_2^\alpha=\gamma_1 \star e_1$. Hence,
its length will be at most the length of the loop $\gamma_1$, which is at most $l+a$ plus the length of the geodesic $e_1$, which is at most $a$.

Note also that the partition $P^\alpha$ can be chosen as fine as desired.

\begin{figure}[! htbp]
\centering
\includegraphics[width=9cm]{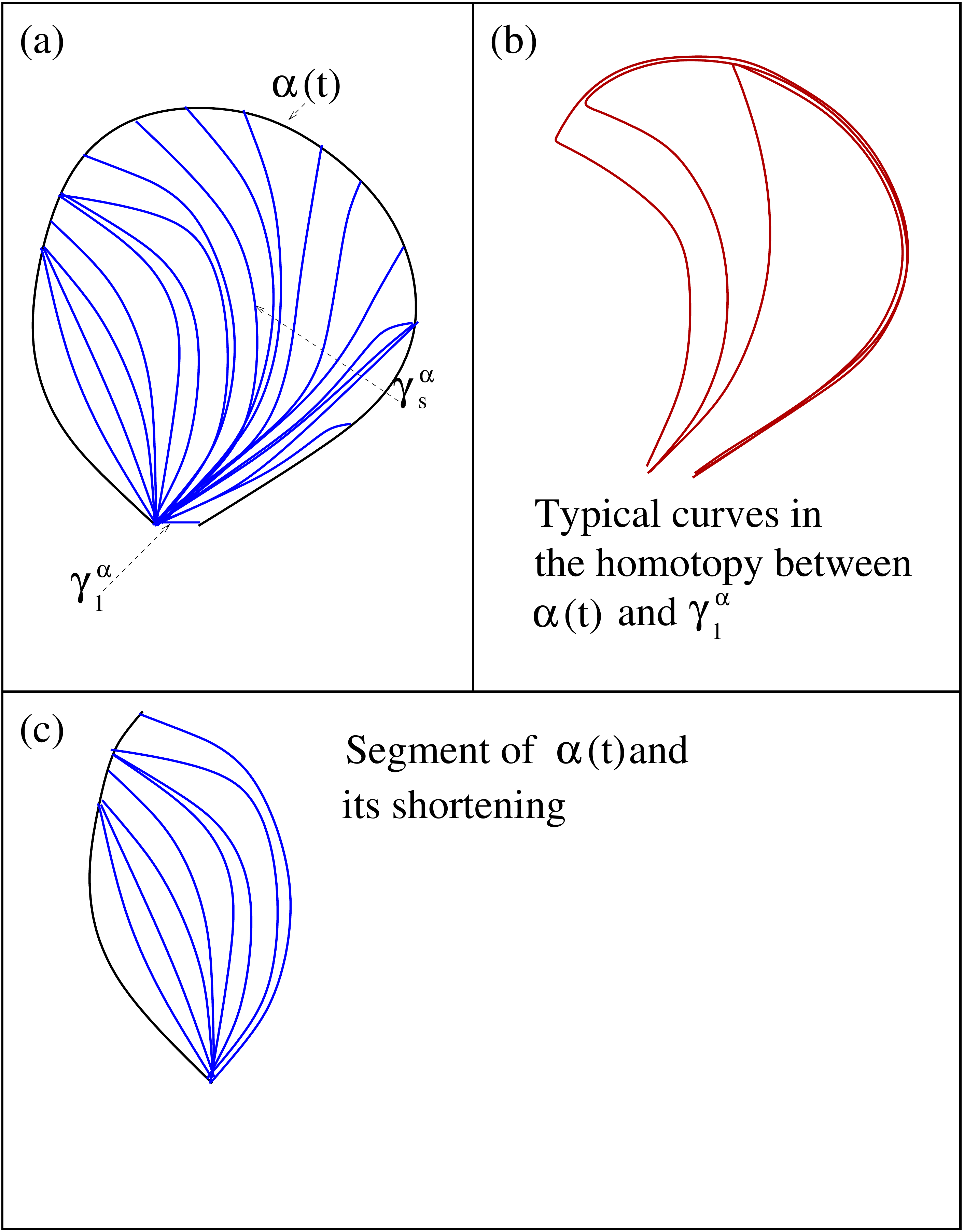}
\caption{Typical curves of the homotopy between $\alpha$ and $\gamma_1^\alpha$}
\label{figure5}
\end{figure}

The one-parameter family $\gamma_s^\alpha$ is schematically pictured  on 
Fig. ~\ref{figure5} (a), while Fig. ~\ref{figure5}(b) depicts the typical curves in the homotopy $H_\alpha(s)$ between $\alpha$ and $\gamma^\alpha_1$. One can see that a typical curve in the homotopy consists of $\gamma^\alpha_s$ followed by the tail part of the 
original curve $\alpha$ that gets shorter and shorter. From now on we will refer to 
these curves as partial shortenings of $\alpha$.
Suppose we would have applied our curve shortening algorithm to a segment of 
$\alpha$, $\alpha|_{[0, x]}$ for some $x \in [0,1]$ as it is depicted on 
Fig. ~\ref{figure5} (c)? It is easy to see that in this case we will obtain 
a $1$-parameter subfamily of the original family $\gamma^\alpha_s$, that is we will get
only the first portion of the curves that connect $p$ with the points of 
$\alpha(t)$ for $t \leq x$.

\medskip

Next we will prove the following theorem,  which  implies 
Theorem ~\ref{TheoremAAA} in the case of $m=1$ by taking $l=2(k-1)a$ and $a=\pi$, which 
in turn implies Theorem ~\ref{Theoremmain} if we take $k={3 \over 2}$.

\begin{Thm} \label{TheoremB}
Let $M^n$ be a closed Riemannian manifold of diameter $d \leq a$ for some positive real number $a$. Let $p \in M^n$ Assume that there exists a positive
number $l$ such that the length of every geodesic loop
that provides a local minimum for the length functional on $\Omega_pM^n$
is {\it not} in the interval $(l, l+2a+\delta]$ for any arbitrarily small $\delta$.
Consider a continuous map
$f:[0,1] \longrightarrow \Omega_p M^n$, (see Fig. ~\ref{figure6} (a)) such that the lengths of
both $f(0)$ and $f(1)$ are at most  $l+a$.
Then $f$ is path homotopic to $\tilde{f}:[0, 1]\longrightarrow
\Omega^{3l+5a+o(1)}_pM^n$.
Moreover, assume that for some $L$ the image
of $f$ is contained in $\Omega_p^LM^n$. Then
one can choose a path homotopy between $f$ and $\tilde f$
so that its image is contained in $\Omega_p^{L+(5a+3l)+o(1)}M^n$.

\end{Thm}

    \begin{Pf}{Proof}

\begin{figure}[! htbp]
\centering
\includegraphics[width=9cm]{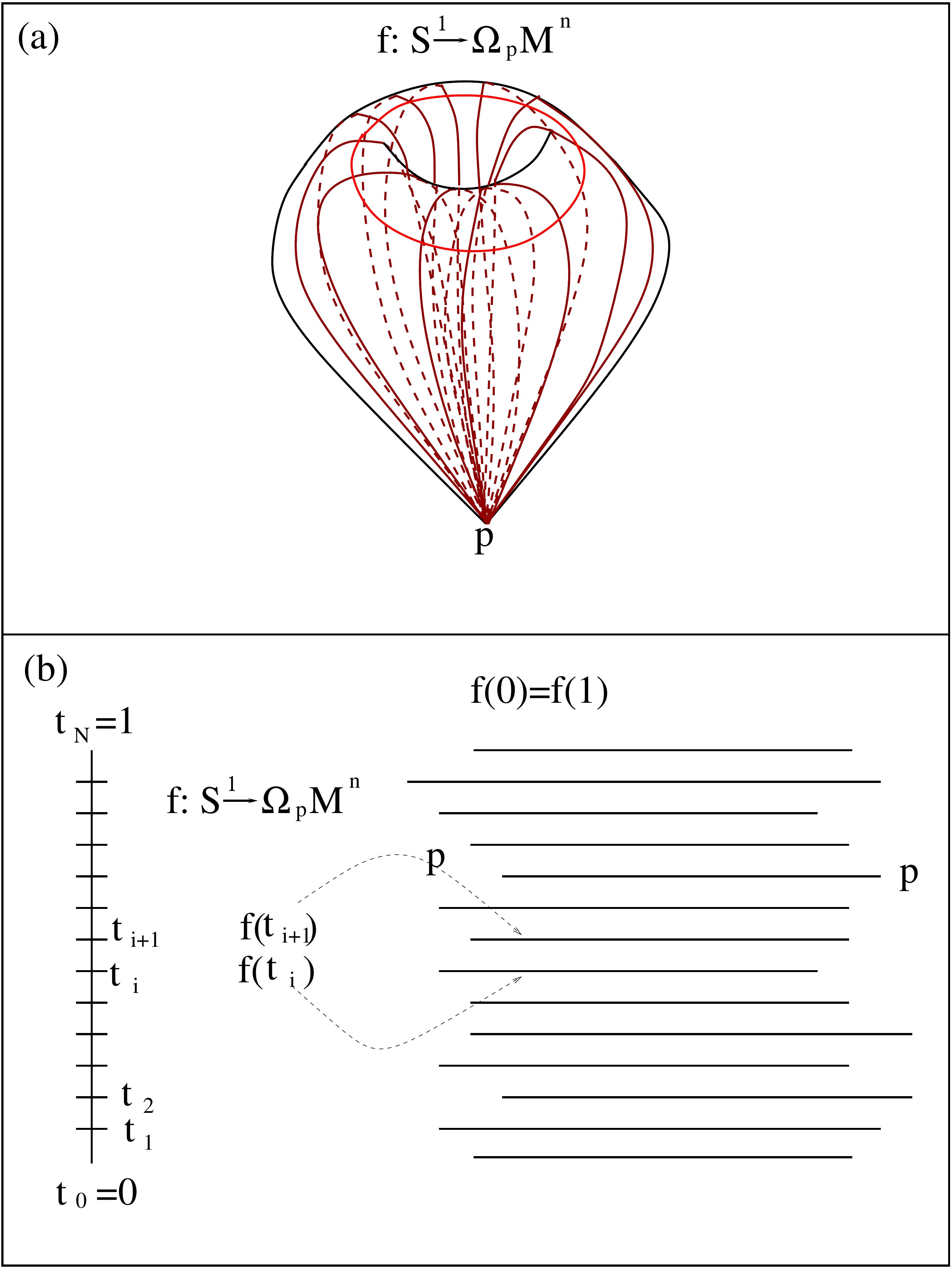}
\caption{Partition of the map $f:S^1 \longrightarrow \Omega_pM^n$ into small intervals}
\label{figure6}
\end{figure}

Without any loss of generality we can assume that $f$ is Lipschitz.
Otherwise, we can make it Lipschitz by an arbitrarily small perturbation.  We will begin by "discretizing" this map in the following sense:
Let   $K=\{t_0=0 < t_1 < t_2 <...<t_N=1\}$ be a partition of 
the interval $[0,1]$, such that 
$\max_i\max_{\tau\in [0,1]}$ length$(f(t)|_{t\in [t_{i-1},t_i]}(\tau)) \leq \epsilon$ for a small 
$\epsilon>0$. We will eventually let $\epsilon$ go to $0$.
Figure ~\ref{figure6} (b) depicts this discretized map $f:S^1 \longrightarrow \Omega_pM^n$, where for simplicity of the picture, we draw $S^1$ and the loops $f(t_i)$ as the segments, 
keeping in mind that (1) $f(0)=f(1)$ and that $f(t_i)(0)=f(t_i)(1)=p$
for all $i=0,...,N$. 

We will now construct a new map $\tilde{f}:S^1 \longrightarrow \Omega_pM^n$ and will show that it is homotopic to the original map $f:S^1 \longrightarrow \Omega_pM^n$. We will do it as follows: first we will replace the loops $f(t_i)$ by the short loops that we obtain via the curve shortening proven in Theorem ~\ref{TheoremA} Let us denote $f(t_i)$ by $\alpha_i$ and 
the "replacement" loops by $\beta_i$. We will "fill in" between $\beta_i$ and 
$\beta_{i+1}$ by contracting $\beta_i*\bar{\beta}_{i+1}$ to $p$ over the loops based at 
$p$ of length at most $2l+4a+\delta+\epsilon$, next we will apply Lemma ~\ref{Lemma} to get a path homotopy between $\beta_i$ and $\beta_{i+1}$ over the loops of length at most 
$3l+5a+2\delta+\epsilon$.
In the construction that follows we will use that the "vertical" curves  $f(t)\vert_{t\in [t_i,\ t_{i+1}]}(\tau)$ that connect ``horizontal" curves $f(t_i)$ are short. 

\begin{figure}[! htbp]
\centering
\includegraphics[width = 9cm]{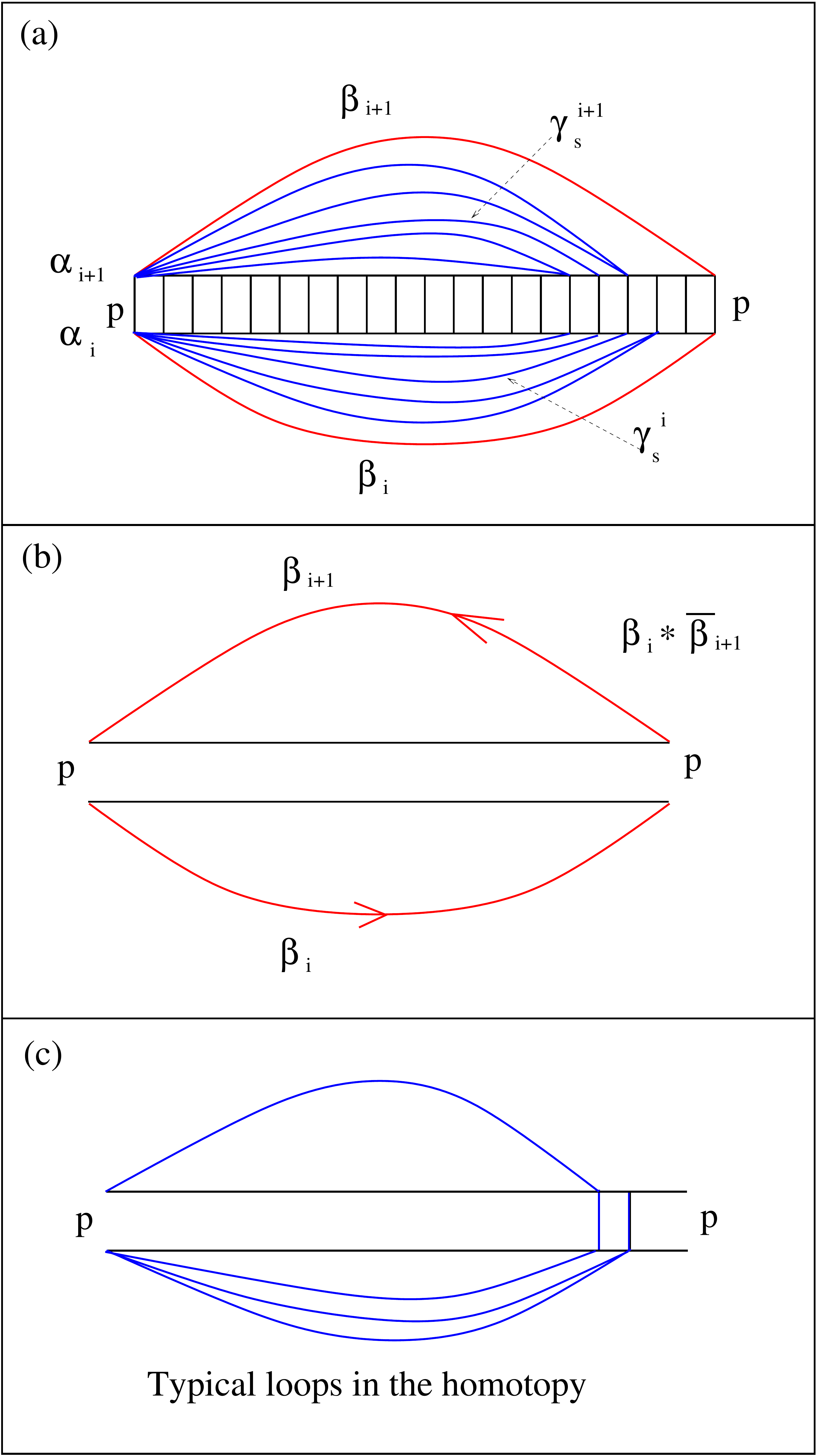}
\caption{Contracting $\beta_i*\bar{\beta}_{i+1}$}
\label{figure7}
\end{figure}

Let us now consider the two loops $\beta_i$, $\beta_{i+1}$ based at the point $p$. 
Despite the fact that $\alpha_i$ and $\alpha_{i+1}$ are loops based at $p$, 
$\beta_i$ and $\beta_{i+1}$ might not be the trivial loops, but their length will be at most $l$, (and not $l+a$ in the case of segments between distinct points $p$ and $q$ ). 
Regardless, of whether or not both $\beta_i$, and $\beta_{i+1}$, or one of the above loops are trivial, we will consider a specific homotopy that is depicted on Fig. ~\ref{figure7} that begins and ends with a constant curve based at $p$. 
Note that by Theorem ~\ref{TheoremA} there exists two $1$-parameter families that will be denoted here as $\gamma_s^i$ and $\gamma_s^{i+1}$
that continuously connect $p$ with the points of $\alpha_i$ and $\alpha_{i+1}$ respectively.
Recall that $\beta_i=\gamma^i_1$, and $\beta_{i+1}=\gamma^{i+1}_1$, (see Fig. ~\ref{figure7}
(a)). We will now consider the loop $\beta_i*\bar{\beta}_{i+1}$ (see Fig. ~\ref{figure7} (b)), (in fact, this curve is a concatenation of two loops based at $p$). Note that 
because the distance between $\alpha_i$ and $\alpha_{i+1}$ is at most $\epsilon$, 
they are continuously connected by the family of curves $f(t)\vert_{t \in [t_{i}, t_{i+1}]}(\tau)$, $\tau \in [0,1]$.

Then there exists a path homotopy between $\beta_i *\bar{\beta}_{i+1}$ and $p$ through the loops 
 $\gamma_{s_1}^i*f(t)\vert_{t\in [t_{i-1}, t_i]}(\tau)*\bar{\gamma}_{s_2}^i$. For this path homotopy   to be defined, the endpoints of $\gamma_{s_1}^i$ and $\gamma_{s_2}^{i+1}$ should be sufficiently 
 $\epsilon$ - close to each other, and to move continuously along $\alpha_i(\tau)$ and
 $\alpha_{i+1} (\tau)$ respectively. Specifically, we want that  
$\gamma_{s_1}^i(1)=f(t_{i})(\tau)$ and 
$\gamma_{s_2}^{i+1}(1)=f(t_{i+1})(\tau)$, or, equivalently,
$\tau_i(s_1)=\tau_{i+1}(s_2)=\tau$. Here $\tau_i(s)$ and $\tau_{i+1}(s)$ denote
the increasing functions from $[0,1]$ to $[0,1]$ for the curves $\alpha_i$ and $\alpha_{i+1}$. Recall that for  every $z$, the inverse image $\tau_j^{-1}(z)$ is either a
point or
a closed interval.  The inverse image  is an interval only when $z$ is a boundary point
for one of the partitions $P^{\alpha_i}$ or $P^{\alpha_{i+1}}$.
Without loss of generality, we may assume that these partitions are disjoint, 
so that the loops
of the form
 $\gamma_{s_1}^i*f(t)\vert_{t\in [t_i, t_{i+1}]}(\tau)*\bar{\gamma}_{s_2}^{i+1}$ form a $1$-parametric family. This family is the desired homotopy.
\par
Note also that, when $\tau_1(s_1)$ is not a boundary point of $P^{\alpha_i}$,
then the length of $\gamma_{s_1}^i \leq l+a+\delta$.
Since $P^{\alpha_i}$ and $P^{\alpha_{i+1}}$ are disjoint, if $\tau_1(s_1)=\tau
=\tau_2(s_2)$, then  at least one of the curves $\gamma_{s_1}^i$ or $\gamma_{s_2}^{i+1}$ is short. That is either the length of $\gamma_{s_1}^i$ is at most $l+a+\delta$
and/or the length of $\gamma_{s_2}^{i+1}\leq l+a+\delta$. Therefore,  each curve
in the homotopy consists of a curve of length at most
$(l+a+\delta)+(l+3a+\delta)+(\epsilon)$. The first summand here is the length of a shorter
curve, the second summand is the length of a potentially longer curve, and $\epsilon$ is the length of the vertical curve connecting the points on the two original loops. Thus the total length of the loops in the homotopy is at most $2l+4a+o(1)$. Here and thereafter, 
$o(1)$ will denote the terms that are bounded by a linear function of $\delta$ and $\epsilon$ that can be made arbitrarily small. The typical loops in this homotopy are depicted on Fig. ~\ref{figure7}(c).

It remains to prove that the constructed path $\tilde f$
in $\Omega_pM^n$ is path homotopic to $f$.

\begin{figure}[!htbp]
\centering
\includegraphics[width=10cm]{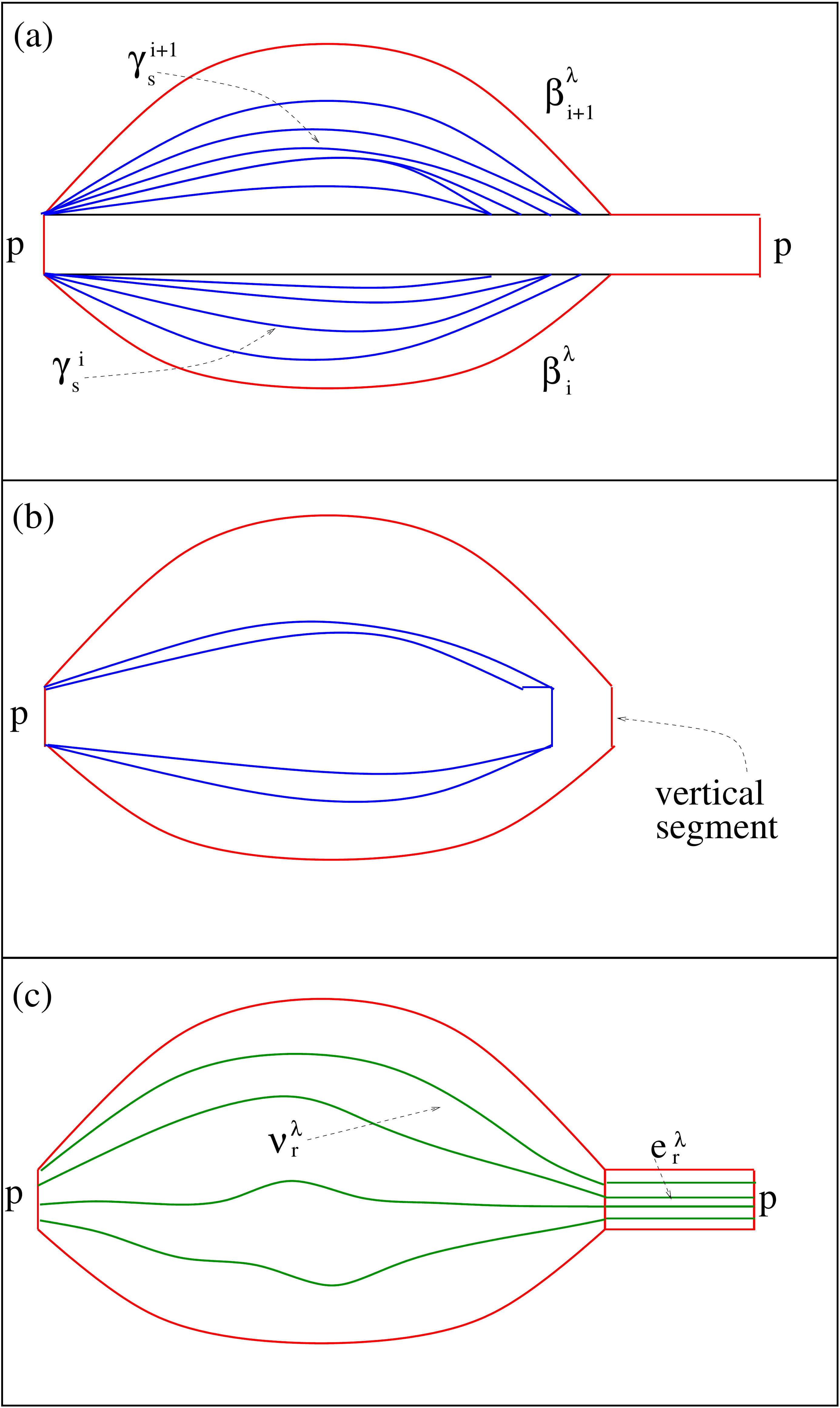}
\caption{Homotopy between $f$ and $\tilde{f}$}
\label{figure8}
\end{figure}

Here is the construction of a path homotopy $G$ between $\tilde f$ and $f$:
$G(1)=\tilde f,\ G(0)=f$. The  main idea to define  $G(\lambda)$  at   $\lambda\in (0,1)$ is to not shorten $f(t_i)$ (for all $i$) all the way
using the construction in the proof of Theorem ~\ref{TheoremA}, as we did
above, but to use the partial shortenings discussed in the proof of 
Theorem ~\ref{TheoremA} instead. We then define $\beta_i^\lambda$, which will consist 
of two arcs. One arc will be $\gamma_{\lambda}^i$ that connects $p$ with $\alpha_i((\tau(\lambda))$
and the second arc is the arc of
$\alpha_i$ that starts at $\alpha_i(\tau(\lambda))$  and ends at 
$\alpha_i(1)$ (see Fig. ~\ref{figure8} (a)).
To construct the desired path homotopy $G$ between $f$ and $\tilde f$
in $\Omega_pM^n$
at the moment $\lambda$ we replace all long
curves $\alpha_i$ not by $\beta_i$ but by $\beta_i^\lambda$.
In other words, $G(\lambda)(t_i)=\beta_i^\lambda$.
Let us consider the two curves $\beta_i^\lambda$ and $\beta_{i+1}^\lambda$ We will construct path homotopy between them that satisfies the following properties: 

\noindent (a) It varies continuously with $\lambda$;

\noindent (b) When $\lambda =0$ we get the original homotopy between 
$f(t_i)$ and $f(t_{i+1})$;

\noindent (c) When $\lambda=1$, we get our homotopy between $\beta_i$ and $\beta_{i+1}$ constructed above.

To construct this homotopy we will assume that both $\alpha_i$ and $\alpha_{i+1}$ are 
parametrized  by a function $\tau(s)$ in such a way that the points $\gamma_s^i(1)$ and $\gamma_s^{i+1}(1)$ are 
$\epsilon$-close to each other. The existence of such a synchronized parametrization
for all the curves in the original homotopy $f(t)$ was discussed in Section 3.4 of [NR7].
We will now join the corresponding points of $\alpha_i$ and $\alpha_{i+1}$ by a short vertical segment, (see Fig. ~\ref{figure8}). 
 $f(t)\vert_{t\in [t_i, t_{i+1}]}(\tau(s))$.
Now we can
form loops
$\gamma_s^i*f(t)\vert_{t\in [t_i, t_{i+1}]}(\tau(s))*\bar{\gamma}_s^{i+1}$
for all $s \leq \lambda$, as before.
These loops will be the loops in a homotopy that contracts the  loop
$\gamma^i_\lambda*f(t)\vert_{t\in [t_{i-1}, t_i]})(\tau(\lambda))*\bar{\gamma}^{i+1}_\lambda$ to a point.

The typical curves of this homotopy are depicted on Fig. ~\ref{figure8}(b).

Let us denote the vertical segment that connects the endpoints of 
$\gamma^i_s$ and $\gamma^{i+1}_s$ as $g_{\tau(s)}(r)$.
We can now construct the curves $G(\lambda)(t)$ for $t \in [t_i,t_{i+1}]$
in the following way. Note that by Lemma ~\ref{Lemma} there exists a path homotopy
between $\gamma^i_\lambda$ and $\gamma^{i+1}_\lambda* g_{\tau(\lambda)}(r)$
over the curves of length at most $6a+3l+o(1)$. This is not quite the homotopy we need. 
Instead, we need the existence of a homotopy between 
$\gamma^i_\lambda$ and $\gamma^{i+1}_\lambda$ over the curves $\nu_r^\lambda$ of 
length at most $6a+3l+o(1)$ such that $\nu^\lambda_r(0)=p$ and $\nu^\lambda_r(1)=g_{\tau(\lambda)}(r)$. The existence proof of such a homotopy is a small modification of the proof of Lemma ~\ref{Lemma}. Thus, it will not be included here. 
We can now extend each such segment $\nu_r^\lambda$ by a segment $e_r^\lambda$, which is a 
tail of the curve from the original homotopy $f(t)$. Thus, when $\lambda=0$, we can see that the curves $f_i$ and $f_{i+1}$ will be unchanged, and so will be the curves interpolating between $f(t)$ and $f(t_{i+1})$. Those, will simply be the curves $f(t)$, whereas, when 
$\lambda=1$, we are replacing the whole curves by $\beta_i, \beta_{i+1}$ respectively, and 
the curves $\nu_r^1$ are the curves $\tilde{f}$. Now, as far as the length bound goes, 
each time we add the tails $e_r^\lambda$ to $\nu_r^\lambda$ we lengthen the curves by at most $L-l-a$ (up to the terms of order $o(1)$). Thus the lengths of curves in the homotopy is bounded by $5a+3l+L+o(1)$.




\par
\par
\par
\par

\par

\end{Pf}

We will next  demonstrate that the above proof generalizes to the maps of spheres of arbitrary dimensions to $\Omega_pM^n$. That is,  if there exist real positive numbers 
$a$ and $l$ such that $d \leq a$ and there exist no geodesic loop based at $p$ of index zero, such that its length is in the interval $(0, l +2a +\delta]$ for some $\delta >0$
then every homotopy class
of $\Omega_pM^n$ can be represented by a sphere that passes
through short loops (thus proving Theorem \ref{TheoremAA}). 
We are going to prove that $f$ is homotopic to a map $\tilde f$
with the image in $\Omega_p^{\tilde L+o(1)}M^n$, where $\tilde L=((4k+2)m+(2k-3))a$, and the homotopy can be chosen so that it has its image in $\Omega_p^{\bar{L}+o(1)}$, where $\bar L=L+\tilde L-(2k-1)a$. We can then take $k={3 \over 2}$ and $a=\pi$

This theorem  is proven  by  recursion with respect to $m$. The case $m=1$ was established in Theorem ~\ref{TheoremB}.
Before summarizing 
the proof of the general case, we will demonstrate it  
in the case when  $m=2$ and $f:S^2 \longrightarrow M^n$. The case of $m=2$ has a straightforward generalization to higher dimensions, yet this dimension is still  small enough as to allow a pictorial description.

\begin{figure}[! htbp]
\centering\includegraphics[width=9cm]{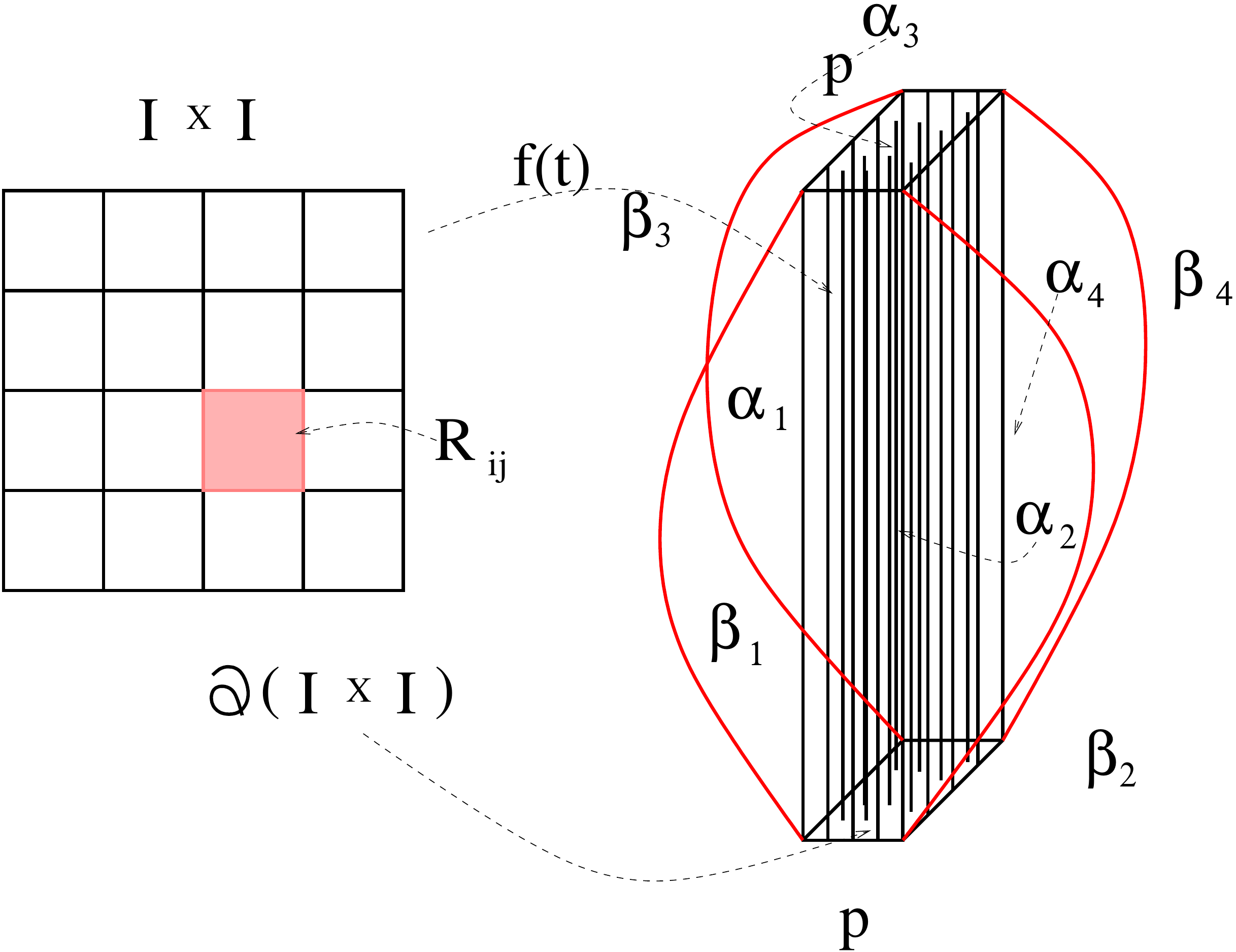}
\caption{Image of a small square in $\Omega_pM^n$}
\label{figure9}
\end{figure}

 Let $I=[0,1]$.  Consider a map
$f:I \times I \longrightarrow \Omega_pM^n$, where $ \partial (I \times I)$
is mapped to $p$.  Without any loss of generality we can assume that all paths in
the image of $f$ are parametrized proportionally to their arclengths. Let us subdivide $I \times I$ into squares of
a very small size. We would like the subdivision be sufficiently small so as  to ensure that the diameter of the image of each small square under  $f$ contributes $o(1)$.    Denote these squares by $R_{ij}$, and their vertices by 
$(x_{i-1},y_{j-1}), (x_{i-1},y_j), (x_i,y_{j-1}), (x_i,y_j)$.  We can consider the subdivision
of $I\times I$ into squares $R_{ij}$ as a cell subdivision. The vertices $(x_i, y_j)$ will be $0$-cells, edges of
squares $R_{ij}$ will be $1$-cells, and their interiors will be $2$-cells. Note that $f$ maps each of
$0$-cells to a loop in $\Omega_pM^n$. 
Replace each of those loops that have length greater than $l+a+\delta$ by a shorter one
as in Theorem ~\ref{TheoremA} via a path homotopy
described in the proof of Theorem ~\ref{TheoremA} (see ~\ref{figure9}).
This yields the desired map $\tilde f$ on the $0$-skeleton of the cell subdivision.

Now we perform the synchronization (as in section 3.4 of [NR7]) of all families $\gamma_s^{f(v)}$, where $v$ runs over the set of all
vertices of all squares $R_{ij}$. This is done in order to find simultaneous parametrizations of all these families by a parameter $s \in [0,1]$ so that
for each $s$ and $v$ the endpoint  $\gamma_s^{f(v)}$ coincides with $f(v)(\tau(s))$, where $\tau(s) \in [0,1]$ does not depend on $v$.

Next,  whenever we apply Lemma 2.1  to a pair of paths that begin at  $p$
and end at some common point naturally associated
with some two vertices of $R_{ij}$, we will refer to those arcs as
$\alpha_1$ and  $\alpha_2$. Note that this choice leads to asymmetry. To achieve the consistency, let us  first 
enumerate all $0$-cells  by numbers $1,2,\ldots$. We will then
denote the path associated with the vertex with a smaller number as $\alpha_1$, and the path associated with the vertex with a larger
number as $\alpha_2$.
\par
Our next step will be to define $\tilde{f}$ on all of the edges that connect these vertices as in the proof of Theorem ~\ref{TheoremB}
in the previous section. This will result in  extending $\tilde f$ to the $1$-skeleton of the cell subdivision to $\Omega_pM^n$. All loops in the image of $\tilde f$ have  the length of at most $(6k-1)a+o(1)$.
What remains now is to define $\tilde{f}$ on the $2$-skeleton of our cell subdivision. Which will be accomplished by 
extending $\tilde{f}$  on the interior of each square $R_{ij}$,  from
$\partial R_{ij}$  for all $i$, $j$.
After this step is accomplished it would remain to verify
that $\tilde f$ and $f$ are homotopic and that this homotopy can be chosen so that all loops in its range are not very long.

\begin{figure}[!htbp]
\centering
\includegraphics[width=9cm]{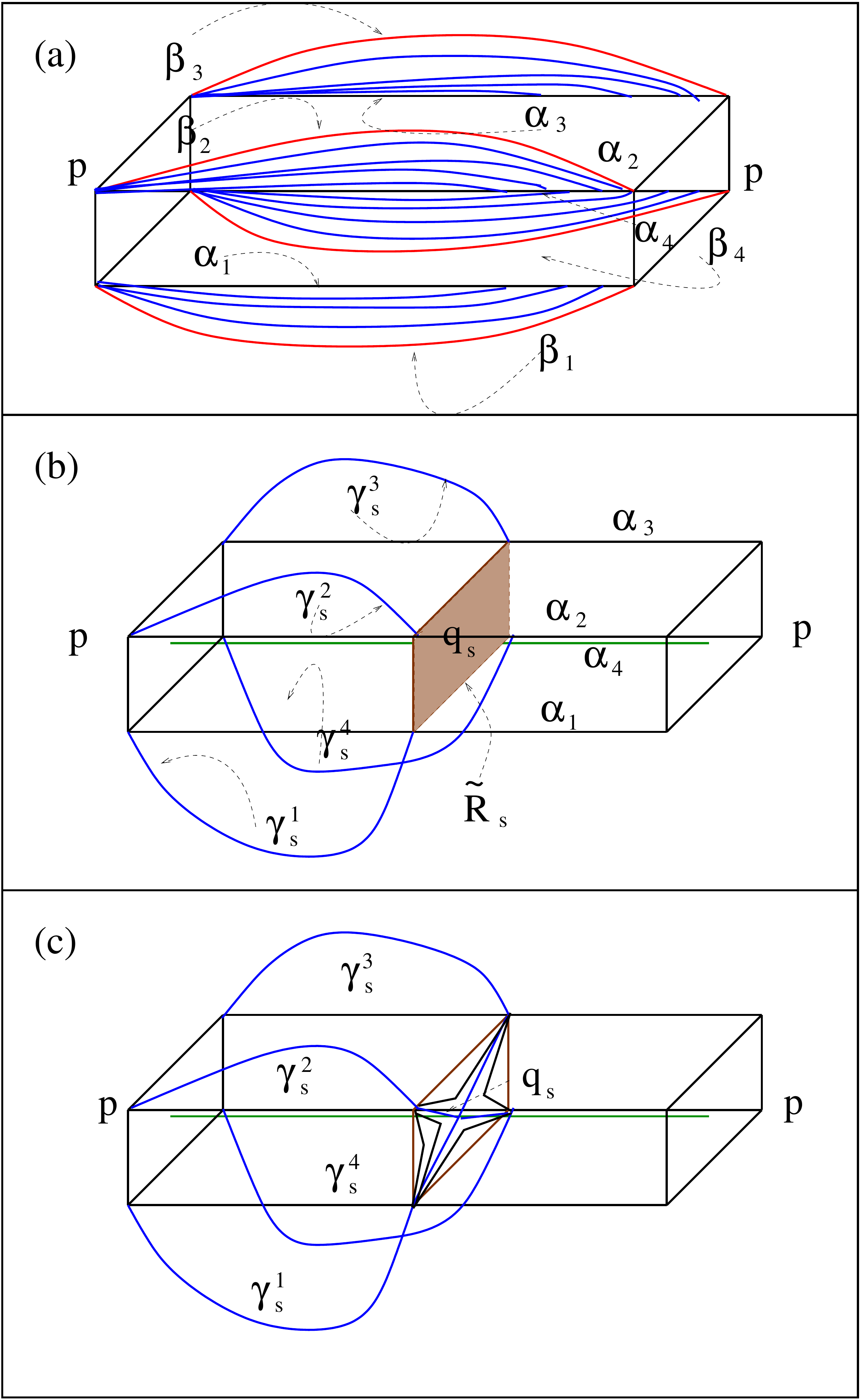}
\caption{Contracting $\tilde{S}^2$ to $p$}
\label{figure11}
\end{figure}

Let us denote the restriction of $\tilde f$ to the boundary 
of a small square $R_{ij}$ by $\tilde{F}$. It induces  a map from $\partial R_{ij} \longrightarrow M^n$. $\partial R_{ij}$ is a $2$-sphere with a natural cell subdivision of the boundary of a parallelipiped.  Two opposite faces of this parallelipiped are mapped to $p$. In our figures below those will be  the left and the right  faces  (see Fig. ~\ref{figure11} (a)).

On Fig. ~\ref{figure11} the two copies of $p$ are depicted as two different points to simplify the diagram. 

Let us consider $R_{ij}$ for some fixed integers $i$ and $j$. Denote 
its vertices  as $v_1, v_2, v_3, v_4$. Let $f(v_k)=\alpha_k$ for 
$k \in \{1,2,3,4\}$, (see Fig. ~\ref{figure11}).
Denote its edges connecting $v_k$ and $v_l$ as $[v_k, v_l]$.  Those edges correspond to the top, bottom, front and back edges of the sphere depicted on  Fig. ~\ref{figure11} (a).   Let us consider how those are mapped to $M^n$ by $\tilde{F}$. 

Without the loss of generality, consider the face that corresponds to the edge $e=[v_1,v_2]$. Let us recall that $\tilde{F}$
was constructed by first replacing $\alpha_1$ and $\alpha_2$ by the curves $\beta_1$, $\beta_2$ respectively. 
We next considered the loop $\beta_1 *\bar{\beta}_2$ based at $p$ and contracted it to  $p$ over the loops of the form 
$\gamma_s^1*f(t)\vert_{t\in e}(\tau(s))*\bar \gamma_s^2$,  $s \in [0,1]$, where
$\tau(s)$ is the synchronization function. Here $\gamma_s^1, \gamma_s^2$ are the curves that connect $p$ with the corresponding points on $\alpha_1(\tau(s)), \alpha_2(\tau(s))$ respectively, the existence of which were proven in Theorem ~\ref{TheoremA}. 

The length of these loops is  at most $2l+4a+o(1)$.
Substituting into this formula $l=2(k-1)a$, we obtain that the length of this loops can be at most $4ka+o(1)$. 

Lemma ~\ref{Lemma} implies the existence of a $1$-parameter family  of short paths interpolating between $\beta_1$ and $\beta_2$ . Similarly we have constructed the map on each 
edge of $\partial R_{ij}$

We will now extend the map to the interior of $R_{ij}$ by first  constructing a $1$-parameter family of $2$-spheres, $\tilde{S}^2_s$ that begins with our sphere and ends with a point. 

This will be done in three steps. 

\noindent {\bf Step 1.} 
For each $s$ consider a $1$-parameter family of "rectangles":$\tilde{R}_s$ that are transversal to 
the curves $f(t)$. One such rectangle is  depicted on Fig. ~\ref{figure11} (b).   Thus, $\tilde{R}_s$ is a $1$-parameter family of rectangles 
that consist of points $f(t)(\tau(s))$ for a fixed $s$ and for all $t \in R_{ij}$.

Next  let us consider  a $1$-parameter family of $2$-spheres 
$S_s^2$, which start with  a $2$ sphere induced by the map 
$f: \partial R_{ij} \longrightarrow \Omega_pM^n$, and end with $S_1^2 =\{p\}$. 
$S_s^2$ then can be described by restricting all of the curves $f(t), t \in \partial R_{ij}$
to the interval $[0,\tau(s),]$. That is, on Fig. ~\ref{figure11} (b), the left part of the 
boundary of the parallelipide, together with the shaded "rectangle" $\tilde{R}_s$
represents $S^2_s$. In order to construct $\tilde{S}^2_s$ we will replace the curves 
by shorter curves as we have done so in the case of $S^2_0$. This construction will continuously depend on the parameter $s$.

Let us partition each $R_s$ continuously with $s$  into families of short curves in the following way. First, let $q_s$ be a $1$-parameter set of centers of these squares. It is represented by a middle line on Fig. ~\ref{figure11} (b). For each $s \in [0,1]$ connect the middle of $R_s$ with each vertex of 
$R_s$. By choosing $R_{ij}$ sufficiently small we can ensure  that 
each connecting edge has length of $o(1)$ and so that the curves very continuously with $s$. Let us denote the edges that we obtain by $\sigma_s^i$ for $i=1,2,3,4$ depending on the order of the vertex. We can now consider the for wedges of curves:
$\sigma_s^1*\bar{\sigma}_s^2, \sigma_s^s*\bar{\sigma}_s^3, \sigma_s^3*\bar{\sigma}_s^4,
\sigma_s^4 * \bar{\sigma}_s^1$. Finally, for each of the wedge and each parameter $s$
we continuously deform it to the corresponding edge via path-homotopy over
the curves of small length, so that they also continuously vary with $s$, (see Fig. 
~\ref{figure11} (c)). 

We will construct $\tilde{S}_s^2$ as follows. First we will replace each segment 
of $\alpha_i$ up to $\tau(s)$ by a corresponding segment $\gamma_s^i$. We will a refer to $\gamma_s^i * \sigma_s^i$ as $\beta_s^i$. That is, we extend each segment $\gamma_s^i$ by a short segment in such a way that the resulting segment connects $p$ with $q_s$. It is possible  
for the same segment to be continuously  replaced by different curves depending on $s$, since $\tau(s)$ is not a one-to-one function. Next, we consider the pairs 
$\gamma_s^j, \gamma_s^{{j+1} \mod 4}$, where $j=1,2,3,4$. In general those curves do not form a loop. However, they will form a loop after  we extend them by  segments 
$\sigma_s^j$. For example, consider $\gamma_s^1, \gamma_s^2$, then the following curve
is a loop based at $p$: $\gamma_s^1 * \sigma_s^1 *\bar{\sigma}_s^2 * \bar{\gamma}_s^2$. 
Moreover, each of these loops can be contracted to $p$ over the loops of length at most 
$4ka+o(1)$, (the first step would be to contract $\sigma_s^1 *\bar{\sigma}_s^2$ to the corresponding edge $f([v_1,v_2])$, while keeping $\gamma_s^1$ and $\bar{\gamma}_s^2$ fixed). 
Now, we can apply Lemma ~\ref{Lemma} to obtain a system of curves interpolating
between $\beta^j_{s}$ and $\beta^{j+1}_{s}$. Those curves will all connect
the point $p$ with the point $q_s$, they will be of length at most $(6k+1)a+o(1)$. 
Thus, we obtain a one parameter family of $2$-spheres in $M^n$ that is endowed with a  sweep-out by short curves connecting $p$ with $q_s$.

\par

\noindent {\bf Step 2.} Now choose a path among these curves and attach it to each of these curves connecting $p$ and $q_s$ so that we will have a family of loops. We will choose this path by first fixing a vertex $v$ of $R_{ij}$, and letting this path be its image under the map that we have constructed in Step 1.   Its length is bounded above by  $ (2k+1)a+o(1)$. We can thus  obtain a sweep-out of each sphere $S^2_s$ by loops of length $\leq (8k+2)a+o(1)$ based at $p$. 
 Finally, define $F$ 
as $F(x,y,t)=f(x,y)(t)$ for all $(x,y)\in R_{ij}, t\in I$.
The image of $v$ under $S^2_t$ is the join of $\gamma^{f(v)}_s$ (from the proof of
Theorem ~\ref{TheoremA}) with the  short straight line segment $\sigma^{F(v)}$
This segment has length of $o(1)$. Vertex $v$ has to be chosen consistently for the sake of continuity with respect to $s$.

To accomplish this consistency we will always choose  $v$ to be the  vertex of $R_{ij}$ with the smallest 
index in the chosen enumeration of all of the  vertices of all of the  squares $R_{ij}$.
As the result we obtain a continuous $1$-parametric family 
of  $\tilde S_t^2$
from $S^1=r_s$ to $\Omega_p^{(8k+2)d+o(1)}M^n$ (parametrized by $t$).
Note that all loops corresponding to a fixed value of $s$ pass through $q_s$.

%
 
\noindent {\bf Step 3.}
Let us  apply Lemma ~\ref{Lemmagen}
to obtain a $3$-disc filling the original $2$-sphere, $S^2_1$. 
At the same time we will obtain a sweep-out of this disk
by loops of length at most 
$(8k+2)a+(2k-1)a+o(1)=(10k+1)a+o(1)$.
In other words, this lemma will produce an extension of $S_1^2$ to a map of a 
$2$-disc to $\Omega_p^{(10k+1)a+o(1)}M^n$ (that can be also regarded as a
map $\tilde F$ of a $3$-disc into $M^n$).

This allows us to extend $\tilde f$ from  $\partial R_{ij}$
to $R_{ij}$. Combining these fillings for all values of $i,j$
we obtain a desired map
$\tilde{f}:S^2 \longrightarrow \Omega_pM^n$  with the image in
$\Omega_p^{(10k+1)a+o(1)}M^n$. 
\par
After the completion of these three steps we will need to prove that the constructed map $\tilde f$ is homotopic to $f$.

\begin{figure}[!htbp]
\centering
\includegraphics[width=9cm]{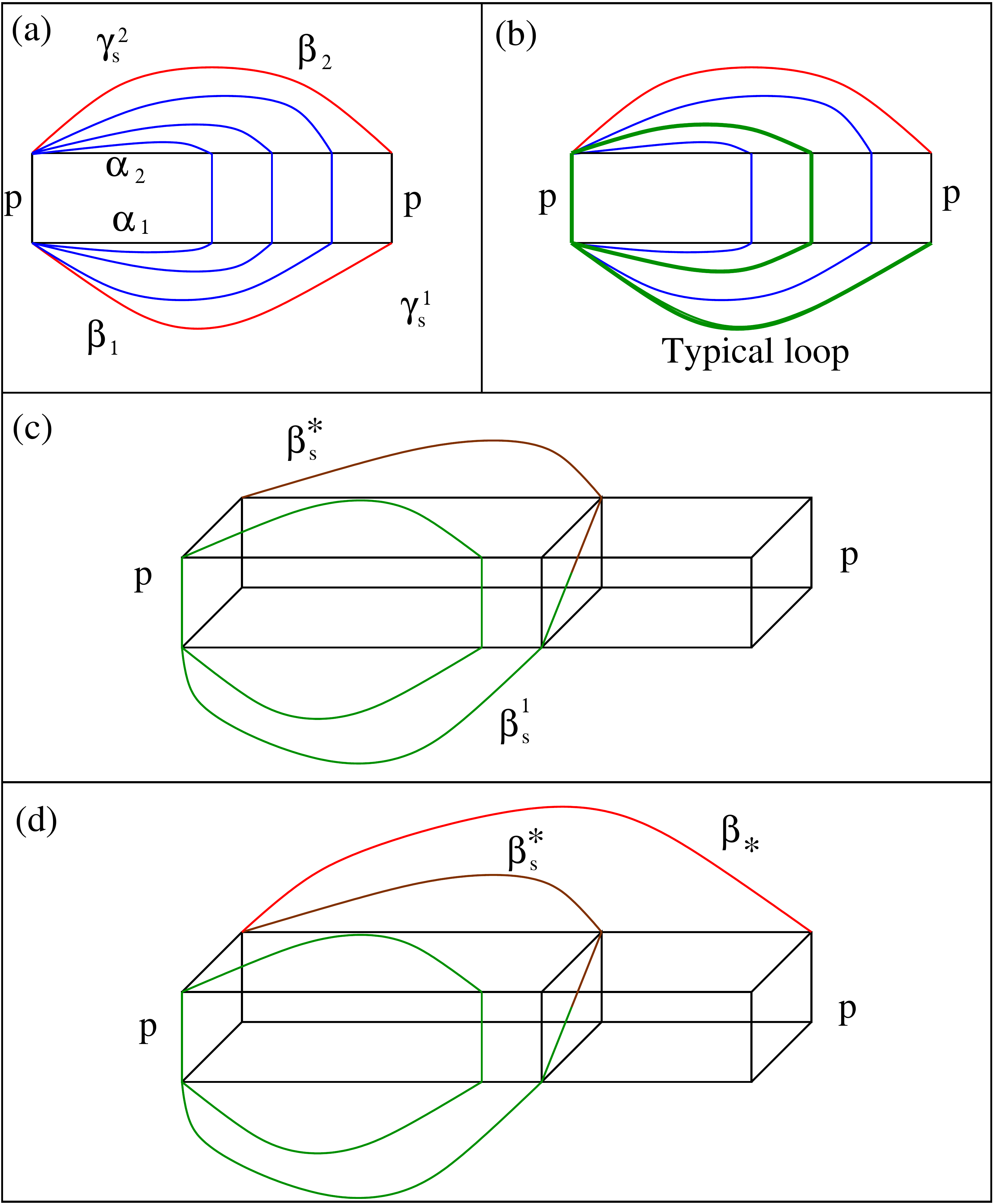}
\caption{Typical Loops}
\label{figure12}
\end{figure}


{\bf A description of two typical loops in the image of $\tilde f$.}
Let us now describe a typical loop in the image of $\tilde{f}$. 
We will consider a loop that corresponds to $t \in R_{ij}$. To understand how
a typical $\tilde{f}(t)$ will "look", we should recall how the map is constructed. 

Recall that the first step  was to change all the curves $\alpha_i=f(v_i)$ for 
$i=1,2,3,4$ to the shorter curves $\beta_i$. Here $v_i$ are the vertices of $R_{ij}$.

At the same step we have constructed families of curves $\gamma_s^i$ that connect
$p$ with the corresponding points on $\alpha_i$. The same point can be connected with 
$p$ by more than one curve (see Fig. ~\ref{figure12} (a)). 

Next we considered a one parameter family of spheres, $S^2_s$. Corresponding to each 
of this spheres, we have constructed a new one parameter family of spheres
$\tilde{S}^2_s$ that consisted of continuously changing with $s$ one-parameter family of "horizontal"  curves that connect the points 
$p$ and $c_s$. Contract the loop $\beta_s^i * \bar{\beta}_s^{(i+1) \mod 4}$ to $p$ and next attach to each such loop  a curve $\beta_s^{v_1}$ that connects $p$ with $c_s$. Those horizontal curves are depicted on Fig. ~\ref{figure12} (c). They are depicted in green and consist of the loop followed by a bottom segment. 
Next we take those curves, and make loops out of them by attaching a segment 
$\beta_s^{v^*}$. This is a top brown curve on Fig. ~\ref{figure12} (c). Thus, the typical loop here would be a loop that consists of a concatenation of loops
$\gamma^i_{\tilde{s}} * vert * \bar{\gamma}^{(i+1) \mod 4}_s$ and $\beta_s^{v^*} * \bar{\beta}_s^{v_1}$. Note that $v^*$ should be chosen canonically, so that the curve changes continuously with respect to $s$.  Next we use Lemma ~\ref{Lemmagen} to construct the loop that fills
the square. It will consist of adding a canonical segment to the loop that we described above (see Fig. ~\ref{figure12} (c)).

\par

%

{\bf Remark: In order to finish the construction in this particular case, we still need to demonstrate that  $f$ and $\tilde f$ are two homotopic maps.}
Denote the desired homotopy by $G$ and its parameter by $\lambda$.
We construct $G$ at $\lambda$ as follows:
Let us first subdivide $I \times I$ into the small rectangles $R_{ij}$ as during the previous
construction.  We can view this as a cell subdivision of $I \times I$, and thus our construction will be inductive with respect to the skeleta of this subdivision. We will 
first consider the vertices and the edges. In this case $G(\lambda)$ will be constructed exactly as in the case of $m=1$.  Recall, that we have done this by changing only the "first" parts of the curves $\alpha_is$ and using their partial shortenings. 
Now it remains to extend $G$ from the boundary of each square $R_{ij}$ to its interior (or, more precisely, from
$\partial R_{ij}\times [0,1]$ to $R_{ij}\times [0,1]$). We are going to assume
that $i,j$ are fixed, and will describe $G(\lambda)\vert_{R_{ij}}:R_{ij}\longrightarrow \Omega_pM^n$ for each value of $\lambda$. 
\par
We will only present an informal idea of this construction. One can consult [NR7] for the details of this proof. 
The informal idea of our construction is that at each moment of time $\lambda$
we would like to define
a map $\tilde F_\lambda$ on $R_{ij}\times [0,s(\lambda)]$ 
the previously defined map $\tilde F$ of $R_{ij}\times I$ restricted to 
$R_{ij} \times [0,s(\lambda)]$, whereas, $\tilde F$ on $R_{ij} \times [s(\lambda), 1]$ should be
the original map $F$ restricted to $R_{ij} \times [s(\lambda), 1]$.
Recall that $\tilde F$ coincides with $\tilde f$ when it is
regarded as a map from $R_{ij}\times I$
into $M^n$, i.e. $\tilde F(x,y,t)=\tilde f(x,y)(t)$.
In particular, we want $\tilde F_1$ to coincide with $\tilde F$.
On the other hand, we would like
the restriction of $\tilde F_\lambda$ on $R_{ij}\times [s(\lambda), 1]$ to coincide with $F$.

The map $G(\lambda)$ will then coincide with $\tilde F_\lambda$, when $\tilde F_\lambda$ is regarded as a
map from $R_{ij}$ to $\Omega_pM^n$.
\par

\begin{figure}[!htbp]
\centering
\includegraphics[width=9cm]{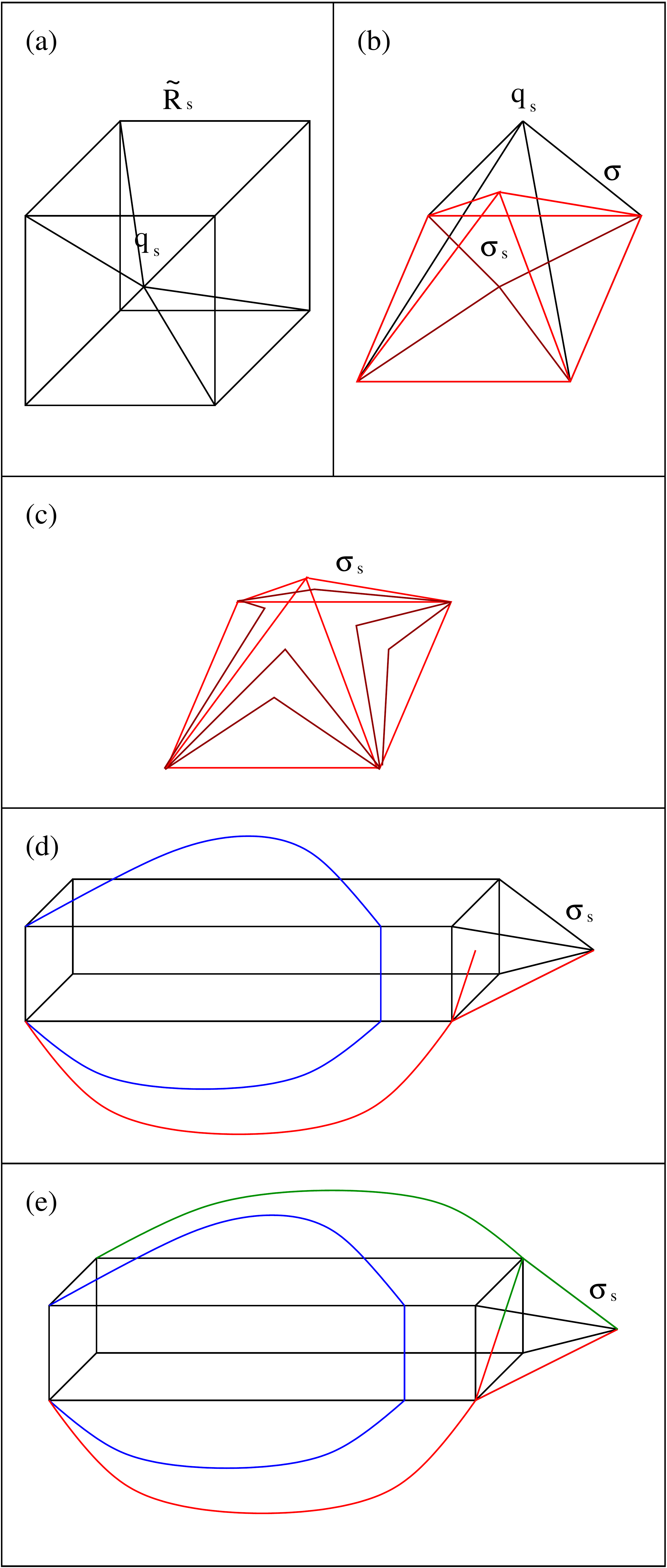}
\caption{General Case}
\label{figure13}
\end{figure}

The lengths of the loops in the image of this homotopy
can be bounded by $L+(10k+3)a+o(1)-(2k-1)a$. We have subtracted 
$(2k-1)a$, because one of the segments of curves $\gamma_s^i$ shortens
a segment whose length is counted as a part of the term $L$

{\bf  Proof of the general case: Construction of $\tilde f$.}
\begin{Pf}{Proof}

\par
We will only present the summary of the proof. One should consult [NR7] for the whole proof, substituting a positive constant $a$ for the diameter $d$

Let $f:I^m \longrightarrow \Omega_pM^n$ be a
continuous map such that $\partial I^m$
is mapped to $p$. 
We will summarize how to construct  $\tilde f:I^m\longrightarrow
\Omega_p^{((4k+2)m+(2k-3))a+o(1)}M^n$.

\par
Let us first consider a map  $F:I^{m+1}\longrightarrow M^n$, induced by $f$.
It is defined as
$f(u)(t)$ for each $u = (x_1,...,x_m) \in I^m$, where $I^m$ corresponds to the first $m$
coordinates of $I^{m+1}$.
Next let us partition $I^m$ into "small" $m$-cubes $R_{i_1,...,i_m}$. 
Let $\tilde{R}_s = F(R_{i_1,...,i_m},s)$ for a fixed $s \in I$
that for each $s$ the image of each straight line in $R_{i_1,\ldots ,i_m}\times\{s\}$ We want the subdivision of $I^m$ be sufficiently to insure that the diameter of 
$\tilde{R}_s$ is smaller than some arbitrarily small $\delta$.  Fig. ~\ref{figure13} (a) depicts such $\tilde{R}_s$. 

The partition of $I^m$ into  small cubes can be viewed as a cell subdivision. Thus, $\tilde{f}$ will be constructed by induction on the order of skeleta. 
By induction assumption, we can construct $\tilde{f}$ on each 
face of  rectangle $R_{i_1,...,i_m}$. In order to extend this map to $R_{i_1,...,i_m}$
we  will construct a $1$-parameter family of $(m-1)$-dimensional spheres, $\tilde{S}_s^{m_1}$ interpolating between  $\tilde{S}_0^{m-1}=\{p\}$ and 
$\tilde{f}(\partial R_{i_1,..,i_m})$. 
So let us fix $s \in I$. Consider  the faces of a $m$-dimensional sphere, which is obtained as follows. We consider $F(R_{i_1,...,i_m},t)$, slice it 
by  rectangle $\tilde{R}_s$, take  the "left" half, and then consider its  boundary 
(see Fig. ~\ref{figure13} (c)). It comes with a natural cellular decomposition, and to construct the new sphere $\tilde{S}^{m-1}_s$ we change the map on each of the cells.  By induction assumption we know how to change the map on 
each of the cells, except for $\tilde{R}_s$. 

In order to extend $\tilde{f}$ to $\tilde{R}_s$, first
partition $\tilde{R}_s$ into  "rectangular pyramides" of dimension 
$m-1$ that has as its base, one of the faces of $\tilde{R}_s$ and as its vertex, the center of 
the cube $q_s$ (see Fig. ~\ref{figure13} (a), (b), (c)). 

Note that the boundary of the pyramide without the base can be continuously deformed to the base, keeping the common sphere  fixed as it is depicted on Fig. ~\ref{figure13} (b). Let us fix the face of the cube. Let $\sigma$ be the corresponding boundary of the pyramide without the base. Consider a deformation of $\sigma$ to the base along $\sigma_r$. On the prior step of the induction, we have obtained a system of 
short curves connecting $p$ and the center of the base, we can continuously repeat this construction by 
deforming the base to $\sigma$ along $\sigma_r$. We can next make the newly constructed segments that correspond to all of the faces of $\tilde{R}_s$ into loops and apply Lemma ~\ref{Lemmagen} (see Fig. ~\ref{figure13} (d), (c)).

\par
\par
\end{Pf}

\par

\noindent {\bf Acknowledgments.} This paper was partially written
during the author's visit of the Max-Planck Institute at Bonn in the Summer 2021.
The author would like to thank the Max-Planck Institute for its
kind hospitality.  The author gratefully
acknowledges partial support by her  NSERC Discovery Grant, RGPIN 220571, during her work on the present paper.

\par

\begin{tabbing}
\hspace*{7.5cm}\=\kill
R. ~Rotman\\
Department of Mathematics\\
40 St. George st.\\
University of Toronto\\
Toronto, Ontario M5S 2E4\\
Canada\\
e-mail: rina@math.toronto.edu

\end{tabbing}
\end{document}